\documentclass[11pt,dvipsnames]{article}
\usepackage{fullpage}
\usepackage[left=1in,top=1in,right=1in,bottom=1in,headheight=3ex,headsep=3ex]{geometry}
\usepackage{graphicx}
\usepackage{pdfpages}
\usepackage[stable]{footmisc}

\usepackage{enumitem}
\usepackage{amssymb}
\usepackage{amsthm}
\usepackage{mdframed}
\usepackage[skip=2pt,font=scriptsize]{caption}

\mdfdefinestyle{E}{
   linecolor=gray!30!white,
  outerlinewidth=0pt,
  roundcorner=20pt,
  innertopmargin=\baselineskip,
  innerbottommargin=\baselineskip,
  innerrightmargin=20pt,
  innerleftmargin=10pt,
  backgroundcolor=gray!30!white}
 \mdfdefinestyle{P}{
    linecolor=gray!5!white,
 frametitlefont=\color{black} \bf ,
 frametitlerulewidth=0pt
  outerlinewidth=0pt,
  innertopmargin=0pt,
  innerbottommargin=0pt,
  innerrightmargin=10pt,
  innerleftmargin=10pt,
  backgroundcolor=gray!3!white}

\title{On the properties of the potential and the spherical harmonic expansions for general planets}
\author{%
O. Costin$^{1}$\quad
R.\,D. Costin$^{1}$\quad
M. Drake$^{1}$\quad
G.\,V. Dunne$^{2}$\quad
C. Ogle$^{1}$\\[0.6em]
\normalsize Corresponding author: O. Costin,
\href{mailto:costin.9@osu.edu}{costin.9@osu.edu}\\[0.8em]
\small $^{1}$Department of Mathematics, The Ohio State University\\
\small 100 Math Tower, 231 West 18th Avenue,
Columbus, OH 43210-1174, USA\\[0.4em]
\small $^{2}$Department of Physics, University of Connecticut\\
\small 196 Auditorium Road, Unit 3046,
Storrs, CT 06269-3046, USA
}
\usepackage[sc]{mathpazo}
\usepackage[T1]{fontenc}
\usepackage{url}
\usepackage{hyperref}\hypersetup{
  colorlinks=true, 
  linktoc=all,  
}
\usepackage{tikz}
\usepackage{amsmath}

\definecolor{labelkey}{rgb}{0,0,1}
\newtheorem{proposition}{Proposition}[section]

\newtheorem{lemma}[proposition]{Lemma}

\newtheorem{remark}[proposition]{Remark}

\mdtheorem[style=P]{theorem}[proposition]{Theorem}
\mdtheorem[style=P]{note}[proposition]{Note}
\mdtheorem[style=P]{definition}[proposition]{Definition}
\mdtheorem[style=P]{definitions}[proposition]{Definitions}
\mdtheorem[style=P]{notations}[proposition]{Notations}
\mdtheorem[style=P]{Proposition}[proposition]{Proposition}
\mdtheorem[style=P]{Note}[proposition]{Note}
\mdtheorem[style=P]{Lemma}[proposition]{Lemma}
\mdtheorem[style=P]{Assumption}[proposition]{Assumption}
\mdtheorem[style=P]{Assumptions}[proposition]{Assumptions}
\mdtheorem[style=P]{Notations}[proposition]{Notations}
\mdtheorem[style=P]{Theorem}[proposition]{Theorem}
\mdtheorem[style=P]{Corollary}[proposition]{Corollary}

\usepackage{hyperref}
\usepackage{array}
\usepackage{xcolor}

\definecolor{clemsonorange}{HTML}{EA6A20}
\hypersetup{colorlinks,breaklinks,linkcolor=clemsonorange,urlcolor=clemsonorange,anchorcolor=clemsonorange,citecolor=black}

\def\RR{\mathbb{R}}

\def\CC{\mathbb{C}}
\def\NN{\mathbb{N}}

\def\DD{\mathbb{D}}

\def\blackslug{\hbox{\hskip 1pt \vrule width 4pt height 8pt depth 1.5pt
\hskip 1pt}}
\def\qed{\quad\blackslug\lower 8.5pt\null\par}

\def\Re{\mathrm{Re\,}}
\def\Im{\mathrm{Im\,}}
\def\CC{\mathbb{C}}
\def\HH{\mathbb{H}}
 \def\RR{\mathbb{R}}
 \def\DD{\mathbb{D}}
 \def\NN{\mathbb{N}}

\def\le{\leqslant}
\def\ge{\geqslant}
\def\epsilon{\varepsilon}
\def\eps{\varepsilon}

\definecolor{red}{rgb}{0.5,0,0}  

 \makeindex
\def\<{\langle}

\def\epsilon{\varepsilon}
\def\phi{\varphi}
 
\def\le{\leqslant}
\def\ge{\geqslant}

\def\>{\rangle}

\def\0{\emptyset}

\usetikzlibrary{arrows.meta,calc}

\begin{document}
\maketitle
\begin{abstract}
The convergence of spherical harmonic expansions (SHEs) of gravitational
potentials below the Brillouin sphere has been a long-standing problem in
potential theory and geodesy. In previous work, some of the present authors
gave a necessary and sufficient condition for convergence along a fixed ray,
formulated in terms of analyticity of a longitudinal average of the balayage
density.

In the present paper, we extend these results to arbitrary physical solid
planets, assuming only that the planet is a compact subset of $\mathbb R^3$
and that its density belongs to $L^\infty$. We establish sharp
regularity properties of the gravitational potential at the Brillouin sphere.
The potential and its first derivative extend continuously to
the sphere, the first derivative is locally H\"older continuous with exponent
$1/2$, and the second derivative satisfies the optimal bound
\[
V''(z)
=
O\bigl(\operatorname{dist}(z,\partial\mathbb D)^{-1/2}\bigr)
\]
under nontangential approach from the exterior.

We also obtain sharp estimates for the spherical harmonic coefficients:
\[
a_\ell=O(\ell^{-3/2}).
\]
If, for almost every latitude, the Brillouin sphere is a one-sided Lebesgue
point from the interior of the longitudinally averaged density, then
\[
a_\ell=o(\ell^{-3/2}).
\]
Both estimates are optimal; in the latter case, the factor multiplying
$\ell^{-3/2}$ can approach zero at an arbitrarily slow rate.

We explicitly recover the balayage density on the Brillouin sphere from the
exterior gravitational potential, and its circular averages from the potential
along a single radial line.

Finally, we give a physical characterization of convergence of the SHE below the
Brillouin sphere. If the SHE converges at such a point, then, after truncating
the planet at any suitable intermediate radius, the gravitational potential
of the removed outer part is completely determined by that of the remaining
truncated planet. This provides a direct interpretation of the strong
constraint imposed by convergence below the Brillouin sphere.
\end{abstract}

\tableofcontents

\section{Introduction}

\paragraph{Background and motivation}

Spherical harmonic expansions (SHEs) of the gravitational potential play a
central role in geodesy, planetary science, and potential theory. For planets,
the classical expansion is known to converge absolutely and uniformly outside
the Brillouin sphere, the smallest sphere centered at the center of mass that
contains the planet.

 What happens below this sphere, however, has remained a subtle and
long-standing problem. The question is also relevant to geophysical modeling,
where downward continuation of the gravitational field relies on the validity
of the SHE at lower altitudes.

\paragraph{Global vs.\ radial SHEs.}
In geodesy, the spherical harmonic expansion of the
gravitational potential is usually expressed as a double series
involving the global spherical harmonic coefficients
$c_{n,m}$ and $s_{n,m}$.
Restricting this expansion to a fixed radial direction
reduces it to the Taylor series at infinity of the potential
along that ray. Throughout this paper, SHE refers to this
one-variable expansion, and its coefficients are called
spherical harmonic coefficients (SHCs).
\paragraph{The convergence problem}

In previous work~\cite{CMP}, some of the present authors resolved this open
problem by providing a necessary and sufficient condition for convergence of
the SHE below the Brillouin sphere along a given ray. Convergence is controlled
by the local regularity of a certain longitudinal average of the balayage
measure near the highest topographic peak along the ray, rather than by global
smoothness of the planet.

More precisely, convergence below the Brillouin sphere occurs if and only if
this longitudinal average is real-analytic at the highest peak. In addition,
they showed that this condition has probability zero in the space of parameters.
In~\cite{CMP}, these results were established under additional assumptions,
including H\"older continuity of the longitudinal average of the balayage
density and the exclusion of polar peaks. In the present
work we show that these assumptions are unnecessary (see
Theorem~\ref{T:iff}); in particular, the required H\"older continuity follows
automatically under our general assumptions (see Note~\ref{sevenp1}).

\paragraph{Purpose and scope of the present work}
In the present paper, we consider the most general class of planets that is
physically meaningful: the planet is an arbitrary compact subset of $\RR^3$
and its density belongs to $L^\infty$. Within this class, we obtain sharp
results for both the regularity of the gravitational potential at the
Brillouin sphere and the decay of its spherical harmonic coefficients.

We first establish sharp regularity properties of the gravitational potential
as the Brillouin sphere is approached. In particular, we obtain optimal bounds
on the growth of its second derivative and show that the potential and its
first derivative extend continuously to the Brillouin sphere, with the first
derivative locally H\"older continuous with exponent $1/2$.

We then show that the spherical harmonic coefficients (SHCs) satisfy
\[
    a_\ell=O(\ell^{-3/2}),
\]
and that this estimate is sharp. 

Under the additional assumption that, for almost every latitude, the
Brillouin sphere is a one-sided Lebesgue point from the interior of the
longitudinally averaged density, the decay improves to
\[a_\ell=o(\ell^{-3/2})\]
This estimate is also sharp.

The improvement reflects the fact that
large-order coefficients are governed by contributions from progressively
thinner layers of the planet near the Brillouin sphere, leading in the limit
to an effective reduction in dimensionality. We also obtain stronger decay
under additional geometric assumptions on the planet.

We next derive explicit reconstruction formulas for the balayage density on
the Brillouin sphere and for its circular averages. In particular, the full
balayage density can be recovered from the exterior gravitational potential,
while its circular averages can be recovered from the potential along a
single radial line. Thus the balayage measure provides a canonical and compact
representation of the exterior gravitational field.

Finally, we obtain a new characterization of convergence below the Brillouin
sphere that reveals a striking nonlocal rigidity. If the SHE converges at a
point below the Brillouin sphere, then, after truncating the planet at any
suitable intermediate radius, the gravitational potential generated by the
removed outer part is completely determined by the potential generated by the
remaining truncated planet. Thus convergence below the Brillouin sphere
imposes a highly nonlocal dependence between the gravitational effects of the
inner and outer portions of the planet. This rigidity provides a direct physical interpretation of the convergence
condition, which requires an exact relation between spatially separated
portions of the planet that cannot reasonably be expected for a realistic
solid celestial body.
\section{Notations and General Assumptions}\label{Sec0}

\paragraph{General notation.}
The unit disk in the complex plane is denoted by
\[
\mathbb{D}:=\{\zeta\in\mathbb{C}:|\zeta|<1\},
\]
and its closure by $\overline{\mathbb D}$. A disk of radius $r$ centered at
$z_0$ is denoted by $\mathbb{D}(z_0,r)$. The open right and left half-planes
are denoted by $\HH^r$ and $\HH^l$, respectively, and $\mathbb B$ denotes the
unit ball in $\RR^3$. In the following,  $\mathbb N=\{1,2,3,\ldots\}$.

For a set $A$, we denote its complement by $A^c$, its interior by $A^\circ$,
its closure by $\overline A$, and its boundary by $\partial A$. The
coefficients $a_\ell$ denote the SHCs of the potential; see
\eqref{eq:eqvs}.

We write $a\lesssim b$ if $a\leq Cb$ for a positive constant $C$
independent of the relevant parameters, and $a\simeq b$ if
both $a\lesssim b$ and $b\lesssim a$.

\begin{Assumptions}\label{A:generic}
We assume that the planet $P$ is a compact subset of $\RR^3$ and that its
density $\rho$ belongs to $L^\infty(P)$. For physical reasons we assume
$\rho\geq0$, although allowing signed densities would cause no difficulty.
We extend $\rho$ by zero outside $P$.\footnotemark
\end{Assumptions}
\footnotetext{Since $P$ is compact, hence Borel measurable, the
extension by zero is measurable and belongs to
$L^\infty(\mathbb R^3)$.}

\begin{Notations}\label{N:notations}
We adopt the convention used in geodesy that the Brillouin sphere, denoted
by $\mathrm{Br}$, is the smallest sphere centered at the center of mass of
$P$ that contains $P$. We place its center at the origin and choose units so
that its radius is $1$.

We denote the maximal radius of the planet at longitude $\lambda$ and
colatitude $\theta$ by
\[
\hat r_m(\theta,\lambda)
   =\sup\{r:(r,\theta,\lambda)\in P\}.
\]

\begin{Note}\label{letitvanish}
\rm
In parts of the paper it will be convenient to set
$\hat r_m(\theta,\lambda)=1$ by formally enlarging $P$ to the Brillouin
sphere and extending the density by zero in the added region.
\end{Note}

Let
\begin{equation}
\label{eq:defrmtheta}
r_m(\theta)=\sup_{\lambda\in[0,2\pi]}\hat r_m(\theta,\lambda).
\end{equation}
Since $P$ is compact, $r_m$ is upper semicontinuous and hence measurable
(see, e.g., \S\ref{semic}). Without loss of generality, we may assume that
$r_m$ does not vanish; otherwise the configuration can be treated as a
collection of disjoint planets. Thus
\[
r_m(\theta)>0.
\]

We define the longitudinally averaged gravitational mass density by
\begin{equation}
\label{eq:deftilderho}
\tilde\rho(r,\theta)
   =G\int_0^{2\pi}\rho(r,\theta,\lambda)\,d\lambda.
\end{equation}
\end{Notations}

\begin{Note}\label{N:Linftytilde}
By Fubini's theorem, $\widetilde\rho$ is measurable,
and Assumption~2.1 implies that
\[
\|\widetilde\rho\|_\infty
\leq 2\pi G\|\rho\|_\infty.
\]
\end{Note}

Densities arising from physical planets have additional properties that may
lead to stronger results. At a point $\mathbf{x}\in\mathbb{R}^3$, the mass
density is classically defined, when the limit exists, by
\[
\rho(\mathbf{x})
=
\lim_{r\to0}
\frac{m(B(\mathbf{x},r))}{\operatorname{Vol}(B(\mathbf{x},r))}.
\]
At the macroscopic scale relevant here, it is natural to model the density of
a physical planet as a bounded function. By the Lebesgue differentiation
theorem, almost every $\mathbf{x}\in P$ is then a Lebesgue point of $\rho$,
that is,
\begin{equation}
\label{defmas}
\lim_{r\to0}
\frac{1}{\operatorname{Vol}(B(\mathbf{x},r))}
\int_{B(\mathbf{x},r)}
|\rho(\mathbf{y})-\rho(\mathbf{x})|\,d\mathbf{y}
=0
\qquad\text{for a.e. }\mathbf{x}\in P.
\end{equation}
For a physical material, density is defined at scales far below the resolution
of gravitational modeling. It is therefore natural to expect the classically
defined density to agree with its local averages at the scales relevant here.

For our results, however, we need only the following weaker condition on the
longitudinally averaged density at the Brillouin sphere.

\begin{Assumption}\label{A:generic1}
In addition to Assumption~\ref{A:generic}, we assume that, for almost every
$\theta\in[0,\pi]$, $1$ is a one-sided Lebesgue point from the interior of
$\tilde{\rho}(\cdot,\theta)$, that is,
\begin{equation}
\label{eq:assumption2}
\lim_{\epsilon\to0^+}
\frac{1}{\epsilon}
\int_0^\epsilon
\left|
\tilde{\rho}(1-s,\theta)-\tilde{\rho}(1,\theta)
\right|\,ds
=0
\qquad\text{for a.e. }\theta\in[0,\pi].
\end{equation}
\end{Assumption}

We occasionally use the following stronger, but still natural, assumption.

\begin{Assumption}\label{A:generic4}
In addition to Assumption~\ref{A:generic}, we assume that, for almost every
$\theta\in(0,\pi)$, $(1,\theta)$ is a one-sided Lebesgue point from the
interior for $\tilde{\rho}$, by which we mean
\[
\frac{1}{2hk}
\int_0^h\int_{-k}^k
\left|
\tilde{\rho}(1-s,\theta-t)-\tilde{\rho}(1,\theta)
\right|\,dt\,ds
\longrightarrow0
\qquad\text{as }h,k\to0^+.
\]
\end{Assumption}
\begin{Note}
Since \eqref{eq:assumption2} is a condition on longitudinal averages of the
density, these averages may have stronger regularity than the density itself,
for instance continuity or absolute continuity. Absolute continuity would
improve the decay rate by an additional factor of $\ell^{-1}$. Whether such
improved regularity occurs for realistic planets could, in principle, be
investigated numerically using high-resolution synthetic models, at degrees
below those at which discretization effects become dominant.
\end{Note}

\paragraph{Remark.}
The distinction between properties holding almost everywhere and everywhere
rarely affects the conclusions of an analytic problem, although there are
notable instances where it does, for example in certain boundary problems for
PDEs. The present setting is another such instance. Both the regularity of the
potential and the asymptotics of its coefficients are sensitive to the behavior
of the density in layers thinning to the Brillouin sphere, whose limiting
surface has zero measure in $\RR^3$. Thus, information that holds almost
everywhere in the volume of the planet need not control the relevant limiting
behavior at the Brillouin sphere.

\subsection{The potential of a planet}

The gravitational potential of a planet $P$ with density $\rho$ is
\begin{equation}
\label{formU}
U(\mathbf z)
=
G\int_P\frac{\rho(\mathbf y)}{|\mathbf z-\mathbf y|}\,d\mathbf y.
\end{equation}

As in~\cite{CMP}, we analyze the potential along a fixed radial line and
choose spherical coordinates so that the line of observation passes from the
center of the planet through the North Pole. The gravitational potential at
the point $\mathbf z=(0,0,z)$ is denoted by
$V(z):=U(0,0,z)$ and is given, for real $z>1$, by
\begin{equation}
\label{eq:potential-planet}
V(z)
=
G\int_P
\frac{\rho(r,\theta,\lambda)\,r^2\sin\theta\,
dr\,d\theta\,d\lambda}
{\sqrt{r^2-2rz\cos\theta+z^2}}.
\end{equation}

The potential $U$ is real analytic in $\RR^3\setminus P$. Moreover,
\eqref{eq:potential-planet}, with the appropriate choice of the square root,
gives an analytic continuation of $V$ to $|z|>1$. In fact, $V$ extends
analytically to a larger domain in $\CC$, as discussed in
\S\ref{2dplanet}.

At infinity, $V$ has the expansion
\begin{equation}
\label{eq:eqvs}
V(z)=\sum_{\ell=1}^\infty a_\ell z^{-\ell},
\end{equation}
which converges absolutely for $|z|>1$; see also~\cite{CMP}.

\section{Main results}

\subsection{Regularity properties of $V(z)$ as $z$ approaches the Brillouin sphere}

In Theorem~\ref{T:BrillouinReg} we show that $V$ and $V'$ extend continuously
to the boundary $|z|=1$, and that $V'$ is locally H\"older continuous with
exponent $1/2$ on $\{|z|\geq 1\}$. The main estimate is
\[
V''(z)
=
O\left(\operatorname{dist}(z,\partial\mathbb D)^{-1/2}\right)
\]
as $z$ approaches the Brillouin sphere nontangentially from the exterior.

\begin{Theorem}\label{T:BrillouinReg}
Consider a planet satisfying Assumption~\ref{A:generic}, with potential $V$
given by \eqref{eq:potential-planet2}.

\begin{enumerate}
\item[(i)]
Let
\begin{equation}
\label{limitsxi}
\xi\in
\left[-\frac{\pi}{2}+a,\frac{\pi}{2}-a\right]
\end{equation}
for some $a>0$, and let $\epsilon_0>0$ be sufficiently small
(for instance, $\epsilon_0<1/12$). Then there exists $C(a)>0$ such that
\begin{equation}
\label{eq:estimv2}
\left|
\epsilon^{1/2}
V''\bigl(e^{i\phi}(1+\epsilon e^{i\xi})\bigr)
\right|
\leq C(a)
\end{equation}
for all $\phi\in(-\pi,\pi]$ and all $\epsilon\in(0,\epsilon_0)$.

In particular, for $\xi=0$ there exists $C>0$ such that
\begin{equation}
\label{eq:estimv20}
\left|
\epsilon^{1/2}
V''\bigl(e^{i\phi}(1+\epsilon)\bigr)
\right|
\leq C
\end{equation}
for all $\phi\in(-\pi,\pi]$ and all $\epsilon\in(0,\epsilon_0)$.

\item[(ii)]
Suppose, in addition, that Assumption~\ref{A:generic4} holds. Then
\begin{equation}
\label{limVsec}
\lim_{\epsilon\to0^+}
\epsilon^{1/2}
V''\bigl(e^{i\phi}(1+\epsilon)\bigr)
=0
\qquad\text{for a.e. }\phi\in(-\pi,\pi].
\end{equation}

\item[(iii)]
The estimates \eqref{eq:estimv2} and \eqref{limVsec} also hold for a
spherical shell when the boundary is approached from inside the shell.

\item[(iv)]
The functions $V$ and $V'$ extend continuously to $\{|z|\geq1\}$.
For an annular planet, they also extend continuously to the inner boundary
of the annulus. Moreover, $V'$ is locally H\"older continuous with exponent
$1/2$ on $\{|z|\geq1\}$.
\end{enumerate}
\end{Theorem}

The following theorem shows that the estimates in
Theorem~\ref{T:BrillouinReg} are optimal under the stated assumptions.

\begin{Theorem}\label{P:decay2}
The exponent $1/2$ in Theorem~\ref{T:BrillouinReg} is optimal: in general,
neither \eqref{eq:estimv20} nor \eqref{limVsec} remains valid if the factor
$\epsilon^{1/2}$ is replaced by $\epsilon^\alpha$ with
$\alpha<1/2$.
\end{Theorem}

The proof of Theorem~\ref{T:BrillouinReg} is given in
\S\ref{T:BrillouinRegProof}, and the proof of
Theorem~\ref{P:decay2} in \S\ref{pf1}.

\subsection{Estimates of the spherical harmonic coefficients}

Precise estimates of the decay of the spherical harmonic coefficients are
important in numerical calculations, in particular for determining truncation
orders when the expansion is evaluated near, or slightly below, the Brillouin
sphere.

The decay of spherical harmonic coefficients of gravitational fields has long
been studied in physical geodesy and approximation theory. Classical potential
theory guarantees convergence of the exterior spherical harmonic expansion
outside the Brillouin sphere, while stronger decay estimates can be obtained
under additional regularity or symmetry assumptions on the mass distribution.
Such assumptions are difficult to justify for a physical celestial body, whose
density may have interfaces, discontinuities, and structure on many spatial
scales.

Here we assume only that the planet is compact and has bounded density. Under
these assumptions we prove the sharp estimate
\[
a_\ell=O(\ell^{-3/2}).
\]
Under Assumption~\ref{A:generic1}, namely that for almost every latitude the
Brillouin sphere is a one-sided Lebesgue point from the interior of the
longitudinally averaged density, this improves to
\[
a_\ell=o(\ell^{-3/2}).
\]
The latter estimate is also sharp: it cannot in general be replaced by
$O(\ell^{-3/2-\epsilon})$ for any $\epsilon>0$, and the factor $o(1)$ may
tend to zero arbitrarily slowly.

\begin{remark}{\rm 

$ $
\begin{itemize}
\item Greater regularity of the topography and density
    can lead to faster decay of the SHCs; see \cite{CMP}.
    \item SHEs are often used slightly below the Brillouin
    sphere, where they generally diverge. In this case,
    truncation at an appropriate order is essential;
    see \cite{BEV}. This is one reason why knowing the decay rate of the SHCs is important.
\end{itemize}}
\end{remark}

The following theorem gives sharp estimates under
Assumption~\ref{A:generic}, and improved estimates under the additional
one-sided Lebesgue-point condition in Assumption~\ref{A:generic1}.

\begin{Theorem}\label{T:opt1}

{\bf A.} Suppose that Assumption~\ref{A:generic} holds.

\begin{enumerate}
\item[(i)]
The SHCs satisfy
\begin{equation}
\label{eq:est0}
a_\ell=O(\ell^{-3/2}).
\end{equation}

\item[(ii)]
The rate in \eqref{eq:est0} is optimal: there exist planets and a subsequence
$\{\ell_n\}_{n\in\NN}$ such that
\[
|a_{\ell_n}|\gtrsim \ell_n^{-3/2}.
\]
\end{enumerate}

{\bf B.} Suppose, in addition, that Assumption~\ref{A:generic1} holds.

\begin{enumerate}
\item[(i)]
Then
\begin{equation}
\label{eq:al}
a_\ell=o(\ell^{-3/2}).
\end{equation}

\item[(ii)]
The rate in \eqref{eq:al} is optimal in the following sense: for any
sequence $\beta_\ell\to0$, with arbitrarily slow decay, there exist densities
satisfying Assumption~\ref{A:generic1}, which may be chosen continuous, and
a subsequence $\{\ell_n\}_{n\in\NN}$ such that
\begin{equation}
\label{eq:lowerbound}
|a_{\ell_n}|
\gtrsim
\beta_{\ell_n}\ell_n^{-3/2}.
\end{equation}
\end{enumerate}

\end{Theorem}

\begin{Note}
A related decay estimate is obtained in
Proposition~\ref{P:decay0}(ii), \S\ref{basedonreg}, under
Assumption~\ref{A:generic4}. That estimate is used in the proof of
Theorem~\ref{P:decay2}.
\end{Note}

The proof of Theorem~\ref{T:opt1} is given in \S\ref{pfTh3p3}.

\subsection{The two-dimensional planet and its complex counterpart}
\label{2dplanet}

Using the zero-extension convention of
Assumption~2.1, we may integrate first in $\lambda$ 
and then use \eqref{eq:deftilderho} and \eqref{eq:deftilderho} to rewrite \eqref{eq:potential-planet} as
\begin{multline}
\label{eq:potential-planet2}
V(z)
=
\int_0^\pi\int_0^{r_m(\theta)}
\frac{\tilde{\rho}(r,\theta)\,r^2\sin\theta}
{\sqrt{r^2-2rz\cos\theta+z^2}}
\,dr\,d\theta
\\
=
\int_T
\frac{\tilde{\rho}(r,\theta)\,r^2\sin\theta}
{\sqrt{r^2-2rz\cos\theta+z^2}}
\,dr\,d\theta
=
\int_T
\frac{\tilde{\rho}(r,\theta)\,r^2\sin\theta}
{\sqrt{z-re^{i\theta}}\sqrt{z-re^{-i\theta}}}
\,dr\,d\theta,
\end{multline}
where the square root of
\[
(z-re^{i\theta})(z-re^{-i\theta})
\]
is taken to be the branch that is positive for large positive $z$. The domain
$T$ in Cartesian coordinates is
\begin{equation}
\label{defcalH}
T
=
\left\{
(r\cos\theta,r\sin\theta):
\theta\in[0,\pi],\ 0\leq r\leq r_m(\theta)
\right\}.
\end{equation}
Its complex counterpart is
\begin{equation}
\label{eq:defT}
T_{\mathbb C}
=
\left\{
re^{i\theta}:
\theta\in[0,\pi],\ 0\leq r\leq r_m(\theta)
\right\}.
\end{equation}
We also introduce the symmetrized set
\begin{equation}
\label{eq:defcalT}
\mathfrak T
=
T_{\mathbb C}\cup T_{\mathbb C}^{*},
\end{equation}
where $*$ denotes complex conjugation, and write
$\mathfrak T^c=\mathbb C\setminus\mathfrak T$.

\begin{figure}[ht]
\centering
\includegraphics[scale=0.3]{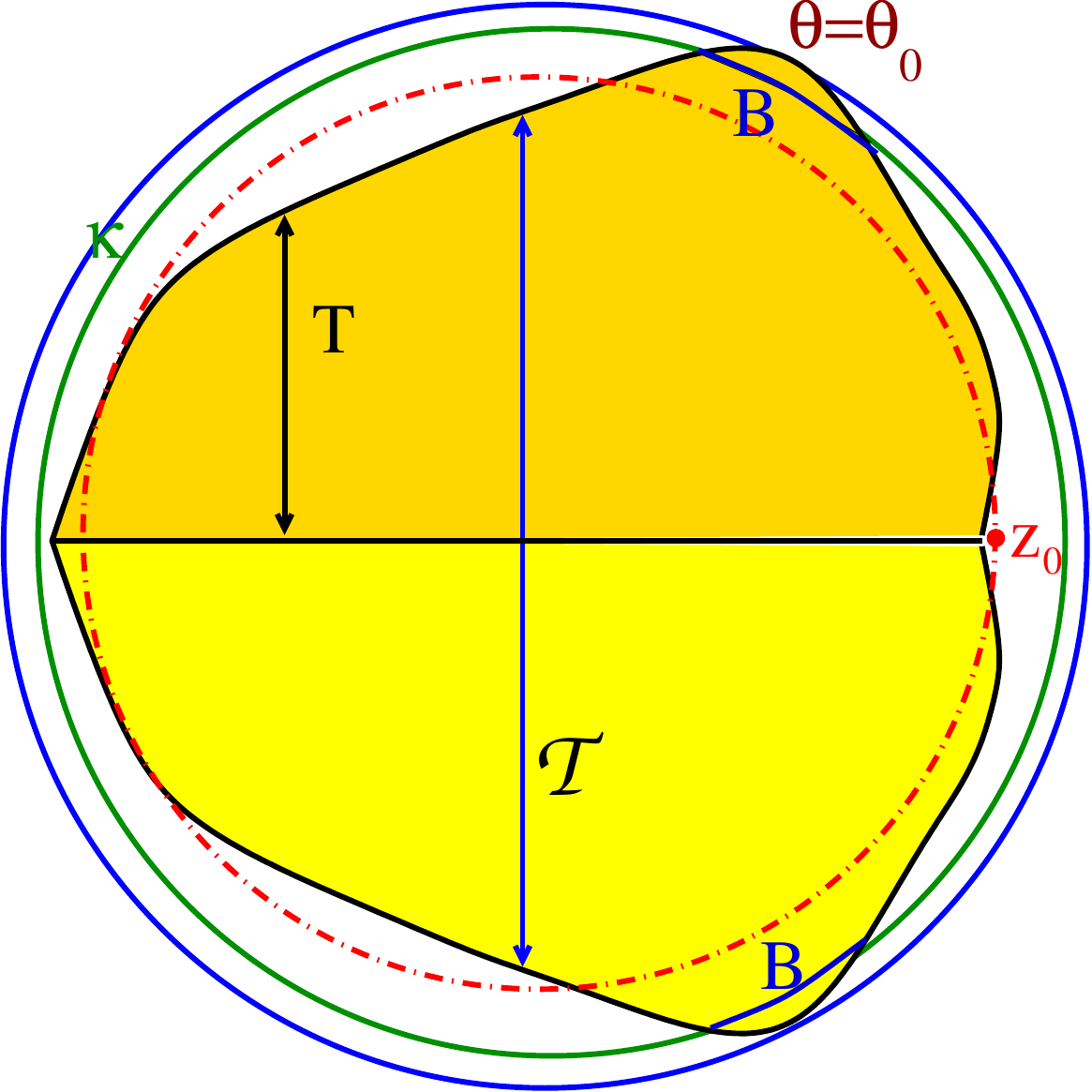}
\caption{The symmetrized two-dimensional planet $\mathfrak T$, consisting
of $T$ (gold) and its complex conjugate $T^*$ (yellow). The Brillouin circle is
shown in blue, and the maximal peak occurs at colatitude $\theta_0$.
Also shown are a point $z_0$ of hypothetical convergence of the SHE and
a circle (green) of radius $\kappa>z_0$.}
\label{fig:flatplanet}
\end{figure}

\begin{Note}\label{N:branch}
By construction, $T$, $T_{\mathbb C}$, and $\mathfrak T$ are star-shaped
with respect to the origin, and hence connected. Their complements are
path-connected.
\end{Note}

The representation \eqref{eq:potential-planet2} shows that $V$ extends
analytically from $|z|>1$ to $\mathfrak T^c$. Indeed, for each
$(r,\theta)\in T$, the two branch points $re^{i\theta}$ and
$re^{-i\theta}$ lie in $\mathfrak T$, and the branch of the square root of
their product is obtained by analytic continuation from infinity.

Thus, after integration in the longitude, the original three-dimensional
problem is reduced to an integral over the two-dimensional domain $T$.
The symmetrized set $\mathfrak T$ describes the obstruction to analytic
continuation of $V$ in the complex $z$-plane; see
Figure~\ref{fig:flatplanet}.

 \begin{remark}{\rm
For axisymmetric planets, the complexified flat planet
$\mathfrak T$ plays a key geometric role: the convergence
properties of the SHE are governed by the complex
singularities of the gravitational potential within
$\mathfrak T$; see \cite{axisymmetric}.
}
\end{remark}

\paragraph{Additional decay in the presence of a sharp maximal peak.}
Suppose that the flat planet $T$ has a unique maximal peak: there is a unique
$\theta_0$ such that
\[
r_m(\theta_0)=1,
\qquad
r_m(\theta)<1 \quad\text{for }\theta\neq\theta_0.
\]
This occurs when the planet $P$ has a unique maximal peak, as is typical for
solid celestial bodies.

If, in addition, the planet satisfies a conical condition near the maximal
peak, then
\[
a_\ell=O(\ell^{-5/2}).
\]
The precise geometric condition and the estimate are given in
Proposition~\ref{P:decay}, \S\ref{pf2}.

\subsection{Can the series converge below Br? Necessary and sufficient conditions}

We now examine the implications of convergence of the SHE below the
Brillouin sphere. Suppose that the series \eqref{eq:eqvs} converges at a
point $z=z_c$ lying above the planet, that is, $z_c\in\mathfrak T^c$, but
below the Brillouin sphere, so that $|z_c|<1$. It then converges for all
$z$ with $|z|>|z_c|$. Theorem~\ref{T:iff} gives necessary and sufficient
conditions for such convergence.

Theorem~\ref{T:iff} reveals a strong rigidity behind convergence below the
Brillouin sphere: after truncating the planet at any suitable intermediate
radius $\kappa$, the potential of the removed outer part is completely
determined by the potential of the remaining truncated planet.

For the following result, assume that the planet is not a sphere and has a
single maximal peak or plateau, centered at colatitude $\theta_0$. We assume
that this maximal feature is isolated: there exists $\delta_0>0$ such that
\[
\sup\{r_m(\theta):|\theta-\theta_0|\geq\delta_0\}<1.
\]
Recall that $V$ is analytic in $\mathfrak T^c$, the exterior of
$\mathfrak T$.

\begin{Notations}
\label{noationa}

Let $z_c\in\mathfrak T^c$, with $|z_c|<1$, and set
\begin{equation}
\label{kappa0}
\kappa_0
=
\max\left(
|z_c|,
\sup\{r_m(\theta):|\theta-\theta_0|\geq\delta_0\}
\right).
\end{equation}

For any $\kappa\in(\kappa_0,1)$, define the truncated planet
\[
P_\kappa
=
\{(r,\theta,\lambda)\in P:r\leq\kappa\},
\]
and let
\[
P^\kappa
=
\overline{P\setminus P_\kappa}
\]
be the closure of the removed outer part.

Similarly, for the flat planet, define
\[
T_\kappa
=
\{re^{i\theta}\in T:r\leq\kappa\},
\qquad
T^\kappa
=
\overline{T\setminus T_\kappa},
\]
and, for the symmetrized domains,
\[
\mathfrak T_\kappa
=
\{z\in\mathfrak T:|z|\leq\kappa\},
\qquad
\mathfrak T^\kappa
=
\overline{\mathfrak T\setminus\mathfrak T_\kappa}.
\]
The choice of $\kappa_0$ ensures that $P^\kappa$ consists only of the outer
part of the maximal peak or plateau and that
$z_c\notin\mathfrak T^\kappa$.

Let $V_\kappa$ and $V^\kappa$ denote the potentials of $P_\kappa$ and
$P^\kappa$, respectively. Equivalently, using
\eqref{eq:potential-planet2}, they are given by the corresponding integrals
over $T_\kappa$ and $T^\kappa$:
\begin{equation}
\label{eq:defdec}
V(z)=V_\kappa(z)+V^\kappa(z),
\end{equation}
where
\[
V_\kappa(z)
=
\int_{T_\kappa}
\frac{\tilde{\rho}(r,\theta)r^2\sin\theta}
{\sqrt{r^2-2rz\cos\theta+z^2}}
\,dr\,d\theta
\]
and
\[
V^\kappa(z)
=
\int_{T^\kappa}
\frac{\tilde{\rho}(r,\theta)r^2\sin\theta}
{\sqrt{r^2-2rz\cos\theta+z^2}}
\,dr\,d\theta.
\]

Denote by
\[
C
=
\partial\mathfrak T^\kappa
\cap
\partial\mathfrak T_\kappa
\]
the common boundary of $\mathfrak T^\kappa$ and $\mathfrak T_\kappa$.
Thus $C$ is the portion of the circle $|z|=\kappa$ separating the truncated
planet from the removed outer part.

Finally, define the annulus
\[
A=\{z:\kappa_0<|z|<1\}.
\]
\end{Notations}

\begin{Note}\label{Note5p2}
The following theorem shows that convergence of the SHE at a point $z_c$
below the Brillouin sphere imposes a strong relation between spatially
separated parts of the planet. For every $\kappa\in(\kappa_0,1)$, the
potential $V^\kappa$ of the removed outer part is completely determined by
the potential $V_\kappa$ of the remaining truncated planet. Consequently,
the potential $V$, and in particular the total mass of $P$, is determined
by $V_\kappa$. Moreover, convergence at $z_c$ requires a continuum of such
conditions, namely \eqref{eq:eqV} for every
$\kappa\in(\kappa_0,1)$. Such a relation cannot reasonably be expected for
a realistic solid celestial body.
\end{Note}
\begin{Theorem}\label{T:iff}
Assume that the SHE \eqref{eq:eqvs} converges at a point $z_c$ with
$0<|z_c|<1$, and use the notation introduced in
Notation~\ref{noationa}.

\begin{enumerate}
\item[(a)]
For every $\kappa\in(\kappa_0,1)$, the potentials $V$ and $V^\kappa$
are completely determined by the potential $V_\kappa$ of the truncated
planet $P_\kappa$.

More precisely, $V_\kappa$ admits an analytic continuation
$\widehat V_\kappa$ to $A\setminus C$, with continuous boundary values on
the two sides of $C$. Denote by $\widehat V_\kappa^+$ the boundary value
obtained by approaching $C$ from the exterior of the circle
$|z|=\kappa$, and by $\widehat V_\kappa^-$ the boundary value obtained
from its interior. Define
\[
\Phi(x)
=
\widehat V_\kappa^-(x)-\widehat V_\kappa^+(x),
\qquad x\in C.
\]
Since $\widehat V_\kappa^+=V_\kappa$ on the exterior side of $C$, this is
equivalently
\[
\Phi(x)
=
\widehat V_\kappa^-(x)-V_\kappa(x).
\]
Then
\begin{equation}
\label{eq:eqV}
V(z)
=
V_\kappa(z)
-
\frac{1}{2\pi i}
\int_C\frac{\Phi(x)}{x-z}\,dx,
\qquad |z|>1,
\end{equation}
where $C$ is oriented counterclockwise. The right-hand side of
\eqref{eq:eqV} gives the analytic continuation of $V$ to $|z|>\kappa$.

\item[(b)]
For every $\kappa\in(\kappa_0,1)$, the total mass of $P$ is completely
determined by $V_\kappa$, and hence by the truncated planet $P_\kappa$.

\item[(c)]
Conversely, if, for some $\kappa\in(\kappa_0,1)$, $V$ satisfies
\eqref{eq:eqV} and its right-hand side is analytic for $|z|>\kappa$,
then the SHE \eqref{eq:eqvs} converges for $|z|>\kappa$.
\end{enumerate}
\end{Theorem}

The proof is given in \S\ref{pf4}.

\subsection{Recovering the balayage measure on the Brillouin sphere}

By the balayage method~\cite{12}, there exists a unique surface mass density
$\sigma$ on the Brillouin sphere $S^2$ whose gravitational potential outside
the sphere agrees with the potential generated by the density $\rho$ on $P$.
We refer to the gravitational field in the region outside the Brillouin sphere as the exterior gravitational field. It has the following integral representation: 
\begin{equation}
\label{Ubalayage1comp}
U(\mathbf z)
=
G\int_{S^2}
\frac{\sigma(\mathbf u)}{|\mathbf z-\mathbf u|}
\,dA(\mathbf u),
\qquad |\mathbf z|>1.
\end{equation}

The balayage density provides a canonical representation of the exterior
gravitational field. Although the inverse gravimetric problem admits
infinitely many mass distributions producing the same exterior potential,
the balayage density is unique once the supporting surface is prescribed.
In this section we derive an explicit relation between $\sigma$ and the
exterior potential $U$. This permits the recovery of $\sigma$ directly from
$U$, whether the latter is obtained from a synthetic model or from
gravitational measurements.

\begin{theorem}\label{T:sigmaU}
Consider a planet satisfying Assumption~\ref{A:generic}, let $U$ be its
potential \eqref{formU}, and let $\sigma\in L^\infty(S^2)$ be its associated
balayage density, defined by \eqref{Ubalayage1comp}. For
$0<\zeta<1$ and $\mathbf u\in S^2$, define
\[
\widetilde U(\zeta,\mathbf u)
:=
G^{-1}\zeta^{-1}U(\zeta^{-1}\mathbf u).
\]
Then the limit
\[
\sigma^*(\mathbf u)
:=
\frac{1}{2\pi}
\lim_{\zeta\uparrow1}
\left[
\zeta\partial_\zeta\widetilde U(\zeta,\mathbf u)
+\frac12\widetilde U(\zeta,\mathbf u)
\right]
\]
exists for every $\mathbf u\in S^2$. Moreover,
\[
\sigma^*=\sigma
\qquad\text{a.e. on }S^2.
\]
\end{theorem}

The proof is given in \S\ref{pfsigmau}.
\subsection{Circular averages}

We next recover the circular averages of the balayage density from the
potential along a single radial line. Circular means on the sphere and their
relation to spherical harmonics are classical; see, e.g., \cite{Muller}.

Define the circular average of $\sigma$ by
\begin{equation}
\label{mudecos}
\mu(\cos\theta)
=
\int_0^{2\pi}\sigma(\theta,\lambda)\,d\lambda.
\end{equation}

The potential $V$ can be expressed in terms of the balayage density $\sigma$
on the unit sphere. Formula~(11) in \cite{CMP} gives\footnote{There is a sign
typo in formula~(11) of \cite{CMP}.}
\begin{equation}
\label{eq:balayage1}
V(z)
=
G\int_{S^2}
\frac{\sigma(\theta,\lambda)\sin\theta\,d\lambda\,d\theta}
{\sqrt{z^2-2z\cos\theta+1}}
=
G\int_0^\pi
\frac{\mu(\cos\theta)\sin\theta\,d\theta}
{\sqrt{z^2-2z\cos\theta+1}}
=
G\int_{-1}^1
\frac{\mu(x)\,dx}
{\sqrt{z^2-2zx+1}},
\qquad z>1.
\end{equation}

It is convenient to move the point at infinity to the origin. Set
$\zeta=1/z$ and define
\begin{equation}
\label{tildeV}
\widetilde V(\zeta)
:=
G^{-1}\zeta^{-1}V(1/\zeta).
\end{equation}
Formula~\eqref{eq:balayage1} expresses $V$ in terms of $\mu$. The following
result gives the converse relation: given $V$, $\mu$ is determined explicitly as follows.

\begin{Theorem}\label{mufromV}
Assume that Assumption~\ref{A:generic} holds; see also
Note~\ref{N:Linftytilde}. Given the potential V, define
\begin{equation}
\label{formH}
H(\zeta)
=
\zeta\widetilde V'(\zeta)
+\frac12\widetilde V(\zeta),
\end{equation}
with $\tilde{V}$ as in \eqref{tildeV}. (The function $H$ is analytic in $\mathbb D$ and extends continuously to
$\overline{\mathbb D}$; cf. Theorem~\ref{T:BrillouinReg}(iv).)

Define $\mu_*:[-1,1]\to\mathbb R$ by
\begin{equation}
\label{muofcost}
\mu_*(\cos\theta)
=
\frac{2}{\pi}
\Im\left[
e^{i\theta}
\int_0^1
\frac{H(te^{i\theta})}
{\sqrt{1-t}\sqrt{1-te^{2i\theta}}}
\,dt
\right],
\qquad 0<\theta<\pi,
\end{equation}
and at the endpoints by
\begin{equation}
\label{muofcostendpoints}
\mu_*(1)=H(1),
\qquad
\mu_*(-1)=H(-1).
\end{equation}
Then $\mu_*$ is continuous on $[-1,1]$ and \eqref{eq:balayage1} is satisfied with $\mu := \mu_*$.
\end{Theorem}

\begin{Note}\label{sevenp1}
\begin{enumerate}
\item[(a)]
If the SHE of $V$ converges at a point $z_c$ below the Brillouin sphere,
then the formula above shows that $\mu_*$ extends analytically to a
neighborhood of $[-1,1]$. This strengthens the corresponding result
in~\cite{CMP}.

\item[(b)]
Formula~\eqref{muofcost} also shows that $\mu_*$ is locally H\"older
continuous with exponent $1/2$ on $(-1,1)$. Thus the H\"older continuity
assumption on the circular average of the balayage density in Theorem~2
of~\cite{CMP} is automatic under Assumption~\ref{A:generic}. Moreover,
the restriction excluding the poles is unnecessary.

\item[(c)] From Assumption \ref{A:generic}, the density $\rho \in L^\infty(P)$, implying the circular average of the balayage density $\mu\in L^\infty([-1,1])$; see
\S\ref{sigmaisLinfty}. In particular, this implies $\mu\in L^2([-1,1])$.
\end{enumerate}
\end{Note}

The proof of Theorem~\ref{mufromV} is given in \S\ref{pfmufromV}.
\section{Proof of Theorem\,\ref{T:BrillouinReg}}\label{T:BrillouinRegProof}

We use the notations in \S\ref{2dplanet}.
For $z\in \mathfrak{T}^c$ we have, by \eqref{eq:potential-planet2}, 
\begin{equation}
\label{dVdz}
\frac{dV}{dz}(z)=-\int_{0}^\pi\, d\theta\int_0^{r_m(\theta)} dr\ \frac{\left( z-r\cos\theta\right)r^2\sin\theta}{\left(  z-re^{i\theta}\right)^{\frac32} \left(  z-re^{-i\theta}\right)^{\frac32} }\, \tilde{\rho}(r,\theta)
\end{equation}

\subsection{Proof of Theorem\,\ref{T:BrillouinReg} (i)}

By allowing the density to vanish (replacing $r_m$ in \eqref{eq:defT} with $1$) we may assume that the planet is the half-disk
\[
T_\mathbb{C}=\DD_+=\{w\in\DD:\Im w\ge 0\}.
\]

\medskip

\noindent\emph{Proof of (i).} 
We have 
\begin{equation}
\label{eq:difdv}
V''(z)=\int_T\tilde{\rho}(r,\theta)
\left(
-\frac{3 r^2 \sin^2\theta}{\left(r^2 - 2 r z \cos\theta + z^2\right)^{5/2}}
+ \frac{2}{\left(r^2 - 2 r z \cos\theta + z^2\right)^{3/2}}
\right)
r^2\sin\theta\, dr\, d\theta.
\end{equation}
Writing $V''(z) = \int_T \tilde{\rho}(r, \theta) h(r, \theta, z) \; dr \; d\theta$, we note the following symmetry: For $z \in \DD^c$, which can be expressed as $z = e^{i\phi}(1+ \epsilon e^{i \xi})$ where $\phi \in (\frac{\pi}{2}, \pi]$, we have

\begin{align}
    V''(e^{i\phi}(1+ \epsilon e^{i \xi})) = \int_T \tilde{\rho}(r, \pi- \theta)h(r,\theta, e^{i(\pi-\phi)}(1+ \epsilon e^{-i \xi})) dr \; d\theta, \label{eq:symv''}
\end{align}
which is $V_*''(e^{i(\pi-\phi)}(1+ \epsilon e^{-i \xi}))$, for the density $\tilde{\rho}_*=\tilde{\rho}(r, \pi- \theta)$.

Note that it suffices to assume $\phi\in[0,\pi/2]$ due to symmetry described in \eqref{eq:symv''}.

Fix a small $\delta>0$ independent of $\epsilon$, and set
\[
D_\delta=\DD_+\cap \DD(e^{i\phi},\delta).
\]
We divide the integration domain $\mathbb{D}_+$ as
\[
\DD_+ = D_\delta \cup (\DD_+\setminus D_\delta).
\]

For $r e^{i\theta} \in \DD_+\setminus D_\delta$, the denominator in \eqref{eq:difdv} is bounded below by a positive constant independent of $\epsilon$. Hence, the integral \eqref{eq:difdv} over this part of the domain is bounded. We therefore restrict our attention to the integral over $D_\delta$, where $r>1-\delta$, and the denominator may become small.

\begin{Lemma}\label{L:ineqsin}

Let $z=e^{i\phi}(1+{\epsilon}e^{i\xi})$ with $\xi$ satisfying \eqref{limitsxi} and $\phi\in (-\pi,\pi]$. Let $\theta \in [0,\pi]$, and $r\in [0,1]$. There exists an $\epsilon_0>0$ small enough so that we have, for all $\epsilon\in(0,\epsilon_0)$
\begin{equation}
\label{eq:eqsint}
\frac{\sin^3\theta}{|z-re^{-i\theta}|^{5/2}}\lesssim 1,
\qquad
\frac{\sin\theta}{|z-re^{-i\theta}|^{3/2}}\lesssim(\sin\theta)^{-1/2},
\end{equation}
and these estimates extend continuously to the limit $\theta\to 0$.
\end{Lemma}

\begin{proof}[Proof of Lemma\,\ref{L:ineqsin}]

 Using
\[
r^2 - 2 r z \cos \theta + z^2=(z-r e^{i\theta})(z-r e^{-i\theta}),
\]
we analyze the case $\phi\in [0,\pi/2]$ (the remaining cases follow by symmetry).

We will prove that
\begin{equation}
    \label{gtsin}
|z-re^{-i\theta}|\gtrsim |\sin\theta|
\end{equation}
which immediately imply inequalities \eqref{eq:eqsint}.

Let $f: [0,1] \to \RR$ be
\[
f({r}) := |z-re^{-i\theta}|^2 = r^2-2r\left[\cos\alpha+\epsilon\cos(\alpha+\xi)\right]+1+2\epsilon\cos\xi+\epsilon^2
,\]
where $\alpha=\theta+\phi\in[0,\frac{3\pi}2]$. 

Let \(w:=1-re^{-i\alpha}\). Since \(\Re w=1-r\cos\alpha\ge 1-r\ge 0.\) Thus, $w$ is in the closed right half plane.
Since \(\xi\in
\left[-\frac{\pi}{2}+a,\frac{\pi}{2}-a\right],\) the angle between \(w\) and \(e^{i\xi}\) (in absolute value) is at most \(\pi-a\). Therefore, \(\Re(we^{-i\xi})
\ge -|w|\cos a.\) (This inequality also holds if $w=0$.) Therefore,
\begin{align}
f(r)
=\left|w+\epsilon e^{i\xi}\right|^2 &=|w|^2+\epsilon^2 +2\epsilon\operatorname{Re}(we^{-i\xi}) \nonumber \\
&\ge |w|^2+\epsilon^2-2\epsilon|w|\cos a \nonumber \\
&=\left(\epsilon-|w|\cos a\right)^2
  +|w|^2\sin^2 a \ge |w|^2\sin^2 a. \label{est1}
\end{align}

It remains only to compare \(|w|\) with \(\sin\theta\). After the symmetry reduction we have \(0\le\phi\le\pi/2\), so \(0\le\alpha=\theta+\phi\le3\pi/2\).

If \(0\le\alpha\le\pi/2\),

$$ |w|^2=(r-\cos\alpha)^2+\sin^2\alpha \ge\sin^2\alpha \ge\sin^2\theta. $$

If \(\pi/2\le\alpha\le3\pi/2\), then \(\cos\alpha\le0\), and hence

$$ |w|\ge\Re w=1-r\cos\alpha\ge1\ge|\sin\theta|. $$

Therefore $f(r)\ge \sin^2 a\sin^2\theta$, which implies

$${ |z-re^{-i\theta}| \ge \sin a\,|\sin\theta|. } $$
\end{proof}

Next, we introduce polar coordinates $(r_2,\theta_2)$ centered at $e^{i\phi}$, with polar axis $te^{i\phi}$, $t\ge 1$. In complex variables,  $D_\delta \subset \{ r_2 e^{i \theta_2}: \theta_2\in(\tfrac\pi2,\tfrac{3\pi}2)\}$, and the coordinate systems are related by
\begin{equation}
\label{eq:iden}
r\sin(\theta-\phi)=r_2\sin\theta_2,\qquad
r\cos(\theta-\phi)=1+r_2\cos\theta_2.
\end{equation}

We rewrite the integral in \eqref{eq:difdv} over $D_\delta$ in terms of $(r_2, \theta_2)$, noting that $r\, dr\,d\theta=r_2\,dr_2\, d\theta_2$. We  compute each term separately. Using Lemma~\ref{L:ineqsin}, and rescaling $r_2=\epsilon\eta$, the first term is bounded by
\begin{align}
\iint_{D_\delta} \tilde{\rho}(r, \theta) \frac{3 r^4 \sin^3\theta}{\left(r^2 - 2 r z \cos\theta + z^2\right)^{5/2}}\; dr \; d\theta &\leq  \|\tilde{\rho}\|_{L^\infty(D_\delta)}\iint_{D_\delta} \frac{3r_2(1+r_2^2+2r_2cos\theta_2)^{3/2}}{|\epsilon e^{i \xi}-r_2e^{i\theta_2}|^{5/2}} \; dr_2 \; d\theta_2 \\
& \lesssim \epsilon^{-1/2}\int_{0}^{\delta/\epsilon}d\eta \int_{ \pi/2}^{3\pi/2}d\theta_2 \, \frac{\eta}{|e^{i \xi} - \eta e^{i \theta_2}|^{5/2}} 
\label{eq:int2}
\end{align}
and noting that $|e^{i \xi} - \eta e^{i \theta_2}|\ge \sin\frac a2\, (1+\eta)$ (since $a\le \theta_2-\xi\le \pi-a$) then in \eqref{eq:int2} we further obtain
\begin{align}
\lesssim \epsilon^{-1/2}\int_0^\infty\frac{\eta\, d\eta}{(\eta+1)^{5/2}}
\lesssim \epsilon^{-1/2}.
\end{align}

For the second term, a similar argument yields
\begin{equation}
\label{eq:int3}
\epsilon^{1/2}\int_0^{\epsilon^{-1}\delta}d\eta\,\int_{\pi/2}^{3\pi/2}d\theta_2\,
\frac{\eta} {|e^{i\xi}+\eta e^{i(\pi-\theta_2)}|^{3/2}}
\lesssim 1
\end{equation}
which completes the proof.

Note that the singular behavior is generated locally near the limiting point: any neighborhood of that point whose radius tends to zero but remains large compared with \(\varepsilon\) captures the full \(\varepsilon^{-1/2}\) scale.

\subsection{Proof of Theorem\,\ref{T:BrillouinReg} (ii) and (iii)}

{ (ii) Since we are looking for an almost everywhere limit,  we do not need to consider the cases $\phi\in\{0,\pi\}$. As in (i) we use symmetries to reduce the analysis to  $\phi\in (0,\pi/2]$. We have already shown that the contribution of the second term of the integrand of \eqref{eq:difdv} to $\epsilon^{1/2}V''((1+\epsilon)e^{i\phi})$ is $O(\epsilon^{1/2})$, so we only need to take the limit of the part of the integral involving the first term:  we evaluate the integral
\[
\int_T \tilde{\rho}(r,\theta)\frac1{(z-re^{i\theta})^{5/2}} g(r,\theta,\epsilon)\, dr\, d\theta
\]
where
\begin{equation}
  \label{eq:defg}
  g(r,\theta,\epsilon)
  =-\frac{3r^4\sin^3\theta}
  {\left(z-r e^{-i\theta}\right)^{5/2}},
  \qquad z=(1+\epsilon)e^{i\phi}.
\end{equation}
 We split the two dimensional planet into the annular sector $(r,\theta)\in D_\delta:=[1-\delta,1]\times [\phi-\delta,\phi+\delta]$ and the rest of the planet, whose contribution is that of an analytic planet with bounds on the derivatives independent of $\epsilon$, and the result is obvious for that part.

We write $(e^{ix}-1)/(ix)=1+h(x)$ where $h(x)$ is analytic and $O(x)$ for small $x$. 

We introduce the scaled variables $r=1-\epsilon\eta,\ \theta=\phi-\epsilon\tau$ and note that
\begin{equation}
  \label{eq:decz}
  z-re^{i\theta}=e^{i\phi}(1+\epsilon-re^{i(\theta-\phi)})=\epsilon e^{i\phi}\left\{1+\eta+ i \tau\left[ 1+h(-\epsilon\tau)\right](1-\epsilon\eta) \right\}
\end{equation}
We obtain
\begin{multline}
  \label{eq:firstt}
  \epsilon^{\frac12}\int_{D_\delta}\int_T \tilde{\rho}(r,\theta)\frac1{(z-re^{i\theta})^{5/2}} g(r,\theta,\epsilon)\, dr\, d\theta
  \\
  =e^{-5i\phi/2}\int_{-\delta/\epsilon}^{\delta/\epsilon} d\tau \int_0^{\delta/\epsilon}d\eta\ \tilde{\rho}(1-\epsilon \eta,\phi-\epsilon\tau)\,K_\epsilon(\eta,\tau) \\
  \text{where }K_\epsilon(\eta,\tau)=\frac{g(1-\epsilon \eta,\phi-\epsilon\tau,\epsilon)}{[1+\eta+\tau i (1+h(-\epsilon\tau))(1-\epsilon\eta)]^{\frac52}}
\end{multline}

Noting that $g(1-\epsilon \eta,\phi-\epsilon\tau,\epsilon)$ has limit as $\epsilon\to 0$ and that for small $\delta$ we have
\[
\left|1+\eta+\tau i (1+h(-\epsilon\tau))(1-\epsilon\eta)\right|\gtrsim (1+\eta+|\tau|  )
\]
it follows
that
\[
|K_\epsilon(\eta,\tau)|
\lesssim
\frac{1}{(1+\eta+|\tau|)^{5/2}} .
\]
on the domain of integration. Extend $K_\epsilon(\eta,\tau)$ by zero to the right half plane $\mathbb H$. Lebesgue's dominated
convergence, therefore, yields
\begin{equation}
\label{eq:KconvL1}
K_\epsilon\longrightarrow K_0(\eta,\tau)=\frac{g(1,\phi,0)}{(1+\eta+i\tau)^{5/2}}
\qquad\text{in }L^1(\mathbb H),
\end{equation}

We use the following elementary consequence of the Lebesgue point
property. If $(1,\phi)$ is a Lebesgue point of $\tilde{\rho}$ and
$K\in L^1(\mathbb H)$, then
\begin{equation}
\label{eq:Lebrescale}
\lim_{\epsilon\to0^+}
\int_{\mathbb H}
K(\eta,\tau)\,
\tilde{\rho}(1-\epsilon\eta,\phi-\epsilon\tau)
\,d\eta\,d\tau
=
\tilde{\rho}(1,\phi)
\int_{\mathbb H}K(\eta,\tau)\,d\eta\,d\tau .
\end{equation}
Indeed, for compactly supported, bounded $K$ this follows directly from
the definition of a Lebesgue point after the change of variables
$s=\epsilon\eta$, $t=\epsilon\tau$.  The general case follows by
approximating $K$ in $L^1$ by bounded compactly supported functions and
using $\tilde{\rho}\in L^\infty$.

Applying \eqref{eq:Lebrescale} to $K_0$, and using
\eqref{eq:KconvL1}, we get
\begin{align}
&\lim_{\epsilon\to0^+}
\int_{\Omega_\epsilon}
\tilde{\rho}(1-\epsilon\eta,\phi-\epsilon\tau)
K_\epsilon(\eta,\tau)\,d\eta\,d\tau
\nonumber\\
&\qquad =
g(1,\phi,0)\,\tilde{\rho}(1,\phi)
\int_0^\infty\int_{-\infty}^{\infty}
\frac{d\tau\,d\eta}
{(1+\eta+i\tau)^{5/2}} .
\label{eq:limitK}
\end{align}
The last integral vanishes. Indeed, it is absolutely convergent, and
for every $\eta\ge0$,
\[
\int_{-\infty}^{\infty}
(1+\eta+i\tau)^{-5/2}\,d\tau
=
-\frac{2}{3i}
\left[
(1+\eta+i\tau)^{-3/2}
\right]_{-\infty}^{\infty}
=0.
\]
Hence the right-hand side of \eqref{eq:limitK} is zero.  Together with
the estimates for the second term and for the complement of
$D_\delta$, this proves
\[
\lim_{\epsilon\to0^+}
\epsilon^{1/2}V''((1+\epsilon)e^{i\phi})=0
\]
for every $\phi$ for which Assumption~\ref{A:generic4} holds, and
therefore for almost every $\phi\in(-\pi,\pi]$.

\

 The proof of {\bf (iii)} is essentially the same, with the exception that some quantities are the complex conjugates of those in case (ii).

\subsection{Proof of Theorem\,\ref{T:BrillouinReg} (iv)}
This follows from the fact that $r \mapsto V''\left( (1+r)e^{i\phi}\right)$ is in $L^1([0,1])$ for every $\phi$ and $V''(z)$ is analytic in $z$ for $|z|>1$.

Indeed, first, $V'$ has radial limits: for $\delta_1>\delta_2$, if $z_{1,2}=e^{i\phi}(1+\delta_{1,2})$
then, using \eqref{eq:estimv2}, we have $\left| V'(z_2)-V'(z_1)\right|=\left|\int_{z_1}^{z_2}V''(w)dw\right|=\left|\int_{\delta_1}^{\delta_2}V''(e^{i\phi}(1+s))ds\right|\le C(a)\int_{\delta_1}^{\delta_2}s^{-1/2}ds\lesssim |\delta_1^{1/2}-\delta_2^{1/2}|$.
Extend $V'$ to $S^1$ by $V'(e^{i \phi}) := \lim_{\delta \to 0}V'(e^{i \phi}(1 + \delta)) $. We still need to prove continuity of $V'$ on $S^1$.

Moreover, we prove that $V'$ is H\"older continuous of exponent $1/2$ on $S^1$. Let $z_{1,2}=e^{i\phi_{1,2}}$ be two points on the circle, close enough, with $|\phi_2-\phi_1|=h\in[0,\epsilon_0]$. Let $w_{1,2}=(1+h) e^{i\phi_{1,2}}$. As before, we write  $V'(z_2)-V'(z_1)=\int_{z_1}^{z_2}V''(w)dw$ where the path of integration is composed of two radial segments, linking $z_j$ to $w_j$, for $j=1,2$ and along the arc circle from $w_1$ to $w_2$. The radial integrals are of order $h^{1/2}$ as we saw before. Also, $|\int_{w_1}^{w_2}V''(w)dw|\le (1+h)\int_{\phi_1}^{\phi_2}|V''((1+h)e^{is})|ds$ which, using again \eqref{eq:estimv2}, is $\le C(a)(1+h)h^{-1/2}(\phi_2-\phi_1)=C_1h^{1/2}$.

Finally, we prove that $V'$ is H\"older continuous of exponent $1/2$ on $\{z;|z|\ge 1\}$. Let $z_{1,2}=e^{i\phi_{1,2}}(1+\delta_{1,2})$ with $|\delta_{1,2}|\le\delta$ small enough. 
Denote $\Delta z=|z_1-z_2|$. If $\min\{\delta_1,\delta_2\}\le2\Delta z$, then
$|\delta_1-\delta_2|\le \Delta z$, hence $\delta_j\le3\Delta z$, and the boundary estimate above,
gives
$
 |V'(z_1)-V'(z_2)|\le C\Delta z^{1/2}
$
If instead $\min\{\delta_1,\delta_2\}>2\Delta z$, the segment that joins
$z_1$ and $z_2$ and stays outside $S^1$ at a distance of at least $\Delta z$.
Therefore, by \eqref{eq:estimv20},
$
 |V'(z_1)-V'(z_2)|
 \le \Delta z\,\sup_{[z_1,z_2]}|V''|
 \le C\Delta z^{1/2}.
$
Thus, $V'$ extends to a locally H\"older continuous function of the exponent
$1/2$ on $\{|z|\ge1\}$.

Finally, $V'$ is uniformly bounded in a neighborhood of $S^1$, therefore integrating $V'$ along radial segments and
circular arcs therefore shows, more simply, that $V$ also extends
continuously to $S^1$.

  \section{Proofs of the estimates for the SHCs}\label{estim SHC}

  We extend the planet to the Brillouin sphere $\mathbb{B}_1$ by setting the density $\rho$ to zero in the space between the topography and its boundary $S^2$, see Note\,\ref{letitvanish}.

 \subsection{Proof of Theorem\,\ref{T:opt1}}\label{pfTh3p3}

 The gravitational potential at $\mathbf{z}=(0,0,z)$ is given by \eqref{eq:potential-planet}.
We now expand the Newtonian kernel $\left(r^2-2rz\cos\theta+z^2\right)^{-1/2}$ as a series in $r/z$ using the generating function of the Legendre polynomials \eqref{eq:gener} and obtain the coefficients in expansion \eqref{eq:eqvs} in the form

\begin{equation}
  \label{eq:SHC10}
  a_{\ell+1}=\int_0^\pi d\theta\, P_{\ell}(\cos\theta) \,\sin\theta\, \int_0^1dr\,\tilde{\rho}(r,\theta)\, r^{\ell+2} 
\end{equation}
where $\tilde{\rho}$ is given by \eqref{eq:deftilderho}.

\begin{proof}
\
{\bf A.}(i) We use \eqref{eq:SHC10}, and \eqref{eq:szegoineq}. We get
\begin{equation}
  \label{eq:SHC1}
 \left| a_{\ell+1}\right|\lesssim \int_0^\pi d\theta |P_{\ell}(\cos\theta)|\int_0^1 r^{\ell+2}dr\lesssim \ell^{-3/2}
\end{equation}
\

{\bf A.}(ii) Consider a spherical planet for which the density is independent of $\lambda$, defined by
\begin{equation}
  \label{eq:defrhho}
 2\pi G\, \rho(1-s,\theta,\lambda)=1+\frac12\sum_{k=2}^\infty\chi_{I_k}(s)\ \ \text{sign}\left[P_{k!}(\cos\theta)\right],\ \ \text{where }I_k=\left[\frac 1{k!},\frac1{(k-1)!}\right)
\end{equation}

We note that the intervals $I_k$ form a disjoint partition of the interval $(0,1)$. Therefore, the series \eqref{eq:defrhho} converges. Also, $\rho\in L^1([0,1],dr)$ and is positive.

We fix $N$ to be large enough and estimate the coefficients $a_{\ell+1}$ in the subsequence $\ell=\ell_n=n!$ for $n>N$. In the following proof, the multiplicative constants absorbed by the symbol "$\lesssim$" are independent of $N$.

In Formula \eqref{eq:SHC10} we split the integral in $r$ into two integrals $\left(\int_{r=0}^{1-N/\ell}dr + \int_{r = 1-N/\ell}^1dr\right)$ that we evaluate separately, yielding the coefficient $a_{\ell+1}$ as a sum of two terms $a_{\ell+1;1}+a_{\ell+1;2}$.

{\em I.} For the term $a_{\ell+1;1}$, noting that $\tilde{\rho} = 2\pi G \rho$, and
\[
0\le \int_0^{1-N/\ell}dr\,\tilde{\rho}(r,\theta)\, r^{\ell+2} \le \int_0^{1-N/\ell}dr\,r^{\ell+2}\left[1+\frac12\sum_{k=2}^\infty\chi_{I_k}(1-r)\right]=\frac32\frac1{\ell+3}\left(1-\frac N\ell\right)^{\ell+3}\lesssim \ell^{-1}e^{-N}
\]
we obtain
\[
\left|a_{\ell+1;1}\right|\lesssim \left(\int_0^\pi d\theta\, \left|P_{\ell}(\cos\theta)\right| \,\sin\theta \right) \ell^{-1}e^{-N} \lesssim \ell^{-3/2}e^{-N}
\]
where for the last inequality we used \eqref{eq:szegoineq}.

{\em II.} The second term, $a_{\ell+1;2}$, is evaluated by defining an  approximation of $\tilde{\rho}$,
\begin{align*}
\tilde{\tilde{\rho}}( \theta) = 1+\frac12\,\text{sign} P_{n!}(\theta)
\end{align*}
 and estimating the error produced, $\tilde{\rho}-\tilde{\tilde{\rho}}$. Note first that, since for $r\in[1-N/\ell,1]$ only intervals $I_k$ with $k\ge n$ appear in integration below (and recalling that $\ell=n!$ and $n>N$) we have 
\[
\int_{1-N/\ell}^1\, r^{\ell+2}\ \left[\tilde{\rho}(r,\theta)-\tilde{\tilde{\rho}}(\theta)\right]\, dr\ \simeq
\int_0^{N/\ell}\, \, (1-s)^{\ell+2}\ \sum_{k=2}^\infty\chi_{I_k}(s)\left[\text{sign} P_{k!}(\theta)-\text{sign} P_{n!}(\theta)\right]\, ds
\]
\[
=\int_0^{N/\ell}\, \, (1-s)^{\ell+2}\ \sum_{k=n+1}^\infty\chi_{I_k}(s)\left[\text{sign} P_{k!}(\theta)-\text{sign} P_{n!}(\theta)\right]\, ds
\]
whose absolute value is less than (see \eqref{eq:defrhho}) $2\int_0^{1/(n+1)!}(1-s)^{\ell+2}=o(\ell^{-1})$. After integration in $\theta$ and using \eqref{eq:szegoineq} we obtain the error term $a_{\ell+1;2}'$ with
\[
\left| a_{\ell+1;2}'\right|= o(\ell^{-3/2})
\]

It remains to estimate  
\[ a_{\ell+1;2}''=-\int_0^\pi d\theta\, P_{\ell}(\cos\theta) \,\sin\theta\,\int_{1-N/\ell}^1dr\,r^{\ell+2}\, \left(1+\frac12\,\text{sign} P_{n!}(\theta) \right)
\]
\[
=\int_0^\pi d\theta\, P_{\ell}(\cos\theta) \,\sin\theta\, \left(1+\frac12\,\text{sign} P_{n!}(\theta) \right)\ \int_{1-N/\ell}^1dr\,r^{\ell+2} 
\]
and, recalling that $\ell=n!$ and $P_\ell$ form an orthogonal set (hence $P_\ell$ is orthogonal to $P_0=1$),
\begin{equation}
    \label{intereq}
 a_{\ell+1;2}''=\int_0^\pi d\theta\, \left|P_{\ell}(\cos\theta)\right| \,\sin\theta\,  \int_{1-N/\ell}^1dr\,r^{\ell+2} =\frac1{\ell+3}\left[1-\left(1-N/\ell\right)^{\ell+3}\right] \int_0^\pi d\theta\, \left|P_{\ell}(\cos\theta)\right| \,\sin\theta
\end{equation}

To evaluate the integral in $\theta$ in \eqref{intereq} we carve out small intervals near the endpoints $0$ and $\pi$ and on the rest of the interval, $J$, we use Hilb's asymptotics \eqref{eq:Piscos}. The contribution at the endpoints is estimated using \eqref{eq:szegoineq},
\[
\int_0^{\ell^{-1}\ln\ell} d\theta\, \left|P_{\ell}(\cos\theta)\right| \,\sin\theta\lesssim  \ell^{-5/2}\ln^{3/2}\ell,
\]
and a similar estimate holds for the integral over $[\pi-\ell^{-1}\ln\ell,\pi]$. On the rest of the interval, $J$, we have
\[
\int_J  d\theta\, \left|P_{\ell}(\cos\theta)\right| \,\sin\theta=\int_J  d\theta\,\sqrt{\frac{2}{\ell\pi\sin\theta}}\ \left|\cos\left[\left(\ell+\frac12\right)\theta-\frac\pi 4\right]+O\left(\ell^{-3/2}\right) \right|\ \sin\theta\gtrsim \ell^{-1/2},
\]
which used in \eqref{intereq} (the contributions from the end points are much smaller) yields
\begin{equation}
    \label{almax}
a_{\ell+1;2}''\gtrsim \ell^{-1/2}\,\frac1{\ell+3}\left[1-\left(1-N/\ell\right)^{\ell+3}\right]\ge \ell^{-1/2}\,\frac1{\ell+3}\left(1-e^{-N}\right)\ge \ell^{-1/2}\,\frac1{\ell+3}\gtrsim \ell^{-3/2}
\end{equation}

Thus, we have $a_{\ell+1}=a_{\ell+1;1}+a_{\ell+1;2}'+a_{\ell+1;2}''$ where $|a_{\ell+1;1}+a_{\ell+1;2}'|\lesssim \ell^{-3/2}e^{-N}$. 
Choosing $N$ large enough it follows that $a_{\ell+1;2}''$ is the largest term, and it satisfies \eqref{almax} and therefore $a_{\ell+1}\gtrsim \ell^{-3/2}\,$ for all $\ell=n!$ with $n>N$.
\end{proof}

{\bf Note.} The above proof yields that
\[
\int_0^\pi |P_\ell(\cos\theta)|\sin\theta\,d\theta \asymp \ell^{-1/2}. \]

\

{\em Proof of }{\bf B.} (i).

We use \eqref{eq:SHC10} for $a_{\ell-2}$, choose a sufficiently large $N$. For $\ell>N+1$ we have
\begin{equation}
  \label{eq:lebesgue0}
  \int_0^1\tilde{\rho}(r,\theta) r^{\ell-1}dr  =\int_0^{1-N/\ell}\tilde{\rho}(r,\theta) r^{\ell-1}dr+\int_{1-N/\ell}^1 \tilde{\rho}(r,\theta) r^{\ell-1}dr
\end{equation}
We look at each integral separately. We have
\begin{equation}
  \label{eq:lebesgue1}
  \left|\int_0^{1-N/\ell}\tilde{\rho}(r,\theta) r^{\ell-1}dr\right|\lesssim \ell^{-1}(1-N/\ell)^{\ell} \lesssim \ell^{-1}e^{-N}
\end{equation}
Using \eqref{eq:szegoineq} and integrating in $\theta$, we see that this term 
contributes with
\begin{equation}
  \label{eq:lebesgue2}
 |a_{\ell-2;1}| \lesssim \ell^{-3/2} e^{-N}
\end{equation}
For the second integral in \eqref{eq:lebesgue0}, changing variable $r=1-\ell^{-1}t$ and integrating by parts we obtain
\begin{multline}
  \label{eq:legcoeff2}
  \int_{1-N/\ell}^1\tilde{\rho}(r,\theta) r^{\ell-1}dr \simeq \ell^{-1}\int_0^Ne^{-t}\tilde{\rho}(1-\ell^{-1}t,\theta)dt=\ell^{-1}\int_0^Ne^{-t}\left(\int_0^t  \tilde{\rho}(1-\ell^{-1}s,\theta)ds\right)'dt\\
  \lesssim \ell^{-1}e^{-N}+\int_0^N dt\  e^{-t}\int_0^t ds\ \tilde{\rho}(1-\ell^{-1}s,\theta)ds
\end{multline}
(where $\int_0^t  \tilde{\rho}(1-\ell^{-1}s,\theta)ds$ is absolutely continuous since $\tilde{\rho}\in L^\infty$). Upon integration in $\theta$, the first term  in \eqref{eq:legcoeff2}  gives rise to a contribution
\begin{equation}
  \label{eq:lebesgue3}
 |a_{\ell-2;2}| \lesssim \ell^{-3/2} N e^{-N}
\end{equation}
To evaluate the second term  in \eqref{eq:legcoeff2} we change the inner variable to $\ell^{-1}s= u$, then we use Assumptions \ref{A:generic} and \ref{A:generic4} and dominated convergence to obtain:
\begin{multline}\label{eq:LebesguePoints}
  \ell^{-1}\int_0^N dt e^{-t}\int_0^t ds \tilde{\rho}(1-\ell^{-1}s,\theta)ds= \ell^{-1}\int_0^N dt e^{-t}\frac{t}{t/\ell} \int_0^{t/\ell} ds \tilde{\rho}(1-u,\theta)du \simeq  \tilde{\rho}(1,\theta)\int_0^N t e^{-t} dt\\=(1-(N+1)e^{-N}) \tilde{\rho}(1,\theta)(1+o(1))
\end{multline}
(if $\tilde{\rho}(1,\theta)\ne 0;$ the case $\tilde{\rho}(1,\theta)=0$ is similar and simpler). Here $o(1)$ is a quantity that goes to zero as $N/\ell\to 0$.  
We now proceed with the integration in $\theta$. We carve out the intervals $[0,\ell^{-1}\log \ell]$ and $[\pi ,\pi -\ell^{-1}\log \ell]$ from the domain
whose total contribution to the integral is, using \eqref{eq:szegoineq},
\begin{equation}\label{eq:a4}
  |a_{\ell-2;3}|\lesssim \ell^{-5/2}  \log^{3/2}\ell 
\end{equation}
On the rest $J$ of the interval we have, using \eqref{eq:Piscos} and the Riemann-Lebesgue Lemma,
\begin{equation}
  \label{eq:asympt}
a_{\ell-2;4}\simeq \ell^{-3/2}\int_{J}d\theta\sin\theta \tilde{\rho}(1,\theta)(1+o(1))\sqrt {\frac {2}{\pi \ell \sin \left(\theta \right)}}\cos \left[\left(\ell +{\tfrac {1}{2}}\right)\theta -{\tfrac {\pi }{4}}\right]+O(\ell^{-5/2})=o(\ell^{-3/2})
\end{equation}
Putting the estimates together, we see that for any $N$ there is an $\ell_N$ such that for all $\ell\ge \ell_N$ we have
we have
$\ell^{3/2} a_{\ell}\lesssim N e^{-N}$. Let $\epsilon>0$ be small enough an choose $N$ so that $Ne^{-N}<\epsilon$. From the estimates above, we conclude that 
\begin{equation}
  \label{eq:finall}
  |\ell^{3/2} a_{\ell}|<\epsilon
\end{equation}
for all large enough $\ell$. 
Since $N$ is arbitrary, then $\epsilon$ is arbitrarily small, and the result follows.

{\em Proof of }{\bf B.} (ii).

  Let $\mu>2$ be a large number. Let $\ell_k$ be any fast growing sequence so that $k\ell_{k-1}/\ell_k\ll\mu^{-k}$ (such as $(k!)!,\ \exp(\exp(\cdots(\exp(k))))$ etc.).

  Consider first the discontinuous $L^{\infty}$ function
  \begin{equation}
    \label{eq:eqlb}
    {\tilde{\rho}}(\theta)=2+\sum_{j\ge 1} \frac1{\mu^j} \,\text{sign}\left[P_{\ell_j}(\cos\theta)\right]
  \end{equation}
   In what follows, the signs $\lesssim$, $\eqsim$ denote inequalities, respectively, equalities, up to positive constants which are independent of $\mu$.

Consider the ball $\mathbb{B}_1$ with axisymmetric, $\lambda$ and $r$-independent density $\tilde{\rho}(\theta)$ (which obviously satisfies Assumptions \ref{A:generic1}). 
From \eqref{eq:SHC10} we see that
\begin{equation}\label{eq:coeffs}
 a_{\ell+1}=\frac1{ \ell+3}\int_0^\pi P_{\ell}(\cos\theta)\sin\theta\, {\tilde{\rho}} (\theta) d\theta 
\end{equation}

We proceed to evaluate the integral in \eqref{eq:coeffs}. 
As in the proof of Theorem \ref{T:opt1} we split the integration interval in \eqref{eq:coeffs} into $I_{k,0}=[0,\sqrt{k}/{\ell_k}]$, $I_{k,\pi}=[\pi-\sqrt{k}/{\ell_k},\pi]$ and $I_k=[0,\pi]\setminus I_{k,0}\setminus I_{k,\pi}$. 

The integral near the end points, on $I_{k,0}$ and $I_{k,\pi}$ is evaluated using \eqref{eq:szegoineq} and the fact that ${\tilde{\rho}}$ is bounded, which
 
\begin{equation}
    \left|\int_ {I_{k,0}}P_{\ell_k}(\cos\theta)\sin\theta\, {\tilde{\rho}} (\theta) d\theta \right|\le \frac{2}{3}\sqrt{\frac{2}{\ell_k \pi}} \, \left(\frac{\sqrt{k}}{\ell_k}\right)^{\frac32}\!\!\|{\tilde{\rho}}\|_\infty \eqsim \frac{k^{3/4}}{\ell_k^2}
\end{equation}
and the integral on $I_{k,\pi}$ is estimated similarly. It follows that $c_0$, the value of the integral over $I_{k,0}\cup I_{k,\pi}$ satisfies
\begin{equation}
    \label{intk0pi}
    |c_0|\lesssim  \frac{k^{3/4}}{\ell_k^2}
    \end{equation}

 To estimate the integral over $I_k$ we use \eqref{eq:eqlb} and obtain
\begin{multline}
\label{eq:serlaTk}
T_k:= \int_{I_k} P_{\ell_k}(\cos\theta)\, \sin\theta\, {\tilde{\rho}} (\theta)\, d\theta =\sum_{j\ge 1} 
\frac1{\mu^j} \, \int_{I_k}  P_{\ell_k}(\cos\theta)\, \sin\theta\  \text{sign}\left[P_{\ell_j}(\cos\theta)\right] d\theta\\
=\sum_{j\ge 1} T_{kj}\ \ \ \text{where}\ \ T_{kj}:=\frac1{\mu^j} \, \int_{I_k}  P_{\ell_k}(\cos\theta)\ \sin\theta\ \text{sign}\left[P_{\ell_j}(\cos\theta)\right]\, d\theta
\end{multline}

We carved out out the endpoints of the interval so that the asymptotic \eqref{Hilb0}, \eqref{HilbPi} applies uniformly on $I_k\cap[0,\pi/2]$, respectively on $I_k\cap[\pi/2,\pi]$.

Denote $c_1=\sum_{j=1}^{k-1} T_{kj}$, $c_2=T_{kk}$ and $c_3=\sum_{j=k+1}^\infty T_{kj}$ which we estimate separately.

 For $j<k$, using \eqref{Hilb0}, \eqref{HilbPi}, we get
\begin{multline}
 \label{Tkj}
T_{kj}\eqsim \frac1{\mu^j\sqrt{\ell_k }} \int_{I_k} {\sqrt { \sin\theta}}\, \cos \left[\left(\ell_k +{\tfrac {1}{2}}\right)\theta -{\tfrac {\pi }{4}}\right]  \, \text{sign}\left[P_{\ell_j}(\cos\theta)\right] d\theta +\mu^{-j}\,{ {o}}\left(\ell_k ^{-1/2}\right)\\
:=\frac1{\mu^j\sqrt{\ell_k }}J_{k,j}+\mu^{-j}\, {o}\left(\ell_k ^{-1/2}\right)
\end{multline}
(we used the fact that integrating in $\theta$ a term of magnitude ${o}\left(\ell_k ^{-1/2}\theta^{-1/2}\right)$ yields a term ${o}\left(\ell_k ^{-1/2}\right)$).

The integrand in $J_{k,j}$ above is piecewise continuous, with piecewise continuous derivative. The points of discontinuity are the zeros of $P_{\ell_j}(\cos\theta)$: denote these zeroes within ${I_k}$ by $t_1,\ldots, t_m$ (in increasing order and with $m\le \ell_j$). We use \eqref{Hilb0}, \eqref{HilbPi}; then we integrate by parts on each interval $[t_{i-1},t_i]$ (with $t_0=\sqrt{k}/{\ell_k},t_{m+1}=\pi-\sqrt{k}/{\ell_k}$):
\begin{multline}
\label{serlajmic}
J_{k,j}=\sum_{i=1}^{m+1}\left\{ {\sqrt { \sin\theta}}\,\frac{\sin\left[\left(\ell_k +{\tfrac {1}{2}}\right)\theta -{\tfrac {\pi }{4}}\right] }{\ell_k +\tfrac {1}{2}} \, \text{sign}\left[P_{\ell_j}(\cos\theta)\right] \right\}\Big|_{\theta=t_{i-1}}^{t_{i}}\\
-\sum_{i=1}^{m+1} \int_{t_{i-1}}^{t_i} \,\frac{\sin\left[\left(\ell_k +{\tfrac {1}{2}}\right)\theta -{\tfrac {\pi }{4}}\right] }{\ell_k +\tfrac {1}{2}} \, \frac1{2{\sqrt { \sin\theta}}}\, \text{sign}\left[P_{\ell_j}(\cos\theta)\right] \,d\theta
 \end{multline}
which implies
\[
\left|J_{k,j}\right|\le \frac{m+1}{\ell_k +\tfrac {1}{2}}+\frac{1}{\ell_k +\frac {1}2}\sum_{i=1}^{m+1} \int_{t_{i-1}}^{t_i}\frac1{2\sqrt{\sin(\theta)}}\,d\theta\le \frac{3(m+1)}{\ell_k +\tfrac {1}{2}}\le \frac{3(m\ell_j+1)}{\ell_k +\tfrac {1}{2}}
\]
therefore
\[
\left|T_{k,j}\right|\lesssim \frac1{\mu^j\sqrt{\ell_k }}\,\frac{3(\ell_j+1)}{\ell_k+\frac12}\ +\mu^{-j}\, {o}\left(\ell_k ^{-1/2}\right)
\]
which implies
\begin{equation}
    \label{intk0pi1}
|c_1|=\left|\sum_{j=1}^{k-1}T_{k,j}\right|\lesssim \ell_k^{-3/2}\sum_{j=1}^{k-1}\frac{\ell_j}{\mu^j}+\sum_{j=1}^{k-1}\frac{1}{\mu^j}{o}\left(\ell_k ^{-1/2}\right)\lesssim \, k\, \ell_{k-1}\ell_k^{-3/2}+ {o}\left(\ell_k ^{-1/2}\right)
\end{equation}

For $j=k$ we have, from \eqref{eq:serlaTk} and using \eqref{eq:szegoineq},
\begin{multline}
    \label{intk0pi2}
c_2:=T_{kk}= \frac1{\mu^k}\int_{I_k} P_{\ell_k}(\cos\theta)\sin\theta\ \text{sign}\,\left(P_{\ell_k}(\cos\theta)\right) d\theta +\mu^{-k}\, {o}\left(\ell_k ^{-1/2}\right)\\
= \frac1{\mu^k}\int_{I_k} \left|P_{\ell_k}(\cos\theta)\right|\sin\theta\, d\theta +\mu^{-k}\, {o}\left(\ell_k ^{-1/2}\right)\eqsim  \frac1{\mu^k} \,\frac1{\sqrt{\ell_k}}+\mu^{-k}\, {o}\left(\ell_k ^{-1/2}\right)=\frac1{\mu^k} \,\frac1{\sqrt{\ell_k}}\left(1+o(1)\right)
\end{multline}
 (As we will see, this term determines the asymptotic behavior of $a_{\ell_k}$).

For $j>k$ we have, from \eqref{eq:serlaTk} then \eqref{eq:szegoineq} that
\[
\left|T_{k,j}\right| \le  \frac1{\mu^j} \,\int_{I_k} \left|P_{\ell_k}(\cos\theta)\right|\sin\theta\, d\theta \eqsim  \frac1{\mu^j} \,\frac1{\sqrt{\ell_k}} 
\]
Therefore
\begin{equation}
  \left| c_3\right|=  \label{intk0pi3}
\left|\sum_{j=k+1}^\infty T_{k,j}\right|\lesssim \ell_k^{-1/2}\sum_{j=k+1}^\infty \frac{1}{\mu^j}\lesssim \, \frac1{\mu^{k+1}}\, \ell_k^{-1/2}
\end{equation}

We obtained that the integral in \eqref{eq:coeffs} is a sum of the four terms $c_0+c_1+c_2+c_3$ which are estimated by \eqref{intk0pi}, \eqref{intk0pi1}, \eqref{intk0pi2} and \eqref{intk0pi3} respectively. We see that $c_0$ and $c_1$ are much smaller than the other two terms. By choosing $\mu$ large enough we see that the term $c_3$ cannot cancel out $c_2$ and the sum behaves like $c_2$.
It follows that
\[| a_{\ell_k}|\gtrsim \mu^{-k} \ell_k^{-3/2}\text{ as } k\to \infty\]
Since $\ell_k$ can be chosen to grow arbitrarily rapidly, the growth of $\mu^k$ is arbitrarily slow relative to  $\ell_k$ and the result follows in the case $\tilde{\rho}\in L^\infty$.

 Replacing the jumps of sign $P_{\ell_k}(\cos\theta) $ by sufficiently steep continuous transitions does not affect the preceding estimates, provided the transition intervals are sufficiently short. The resulting uniformly convergent series defines a continuous density.
  We omit the straightforward details.

 We showed that the $\beta_{\ell_k}$ in \eqref{eq:lowerbound} satisfies $\beta_{\ell_k}\eqsim \mu^{-k}$. The sequence $\{\ell_k\}_k$ in the proof can be replaced by any faster growing sequence. This means that $\beta_{\ell_k}$ can decrease arbitrarily slowly to zero.

\subsection{Proof of Theorem \ref{P:decay2}}\label{pf1}

\subsubsection{Estimates of the coefficients of the series \eqref{eq:eqvs} based on the regularity of $V$}\label{basedonreg}

The key idea of the proof of Theorem \ref{P:decay2} is to show that the estimates for the decay of the SHC's can be obtained from the regularity of $V$.

\begin{Proposition}\label{P:decay0} 

 (i) Under Assumption~\ref{A:generic}, the coefficients of the series \eqref{eq:eqvs} satisfy
\begin{equation}
\label{eq:coefcnt00}
a_{\ell}=O(\ell^{-3/2}).
\end{equation}

(ii) Under Assumptions \ref{A:generic4}, the coefficients of the series \eqref{eq:eqvs} satisfy
\begin{equation}
\label{eq:coefcnt002}
a_{\ell}=o(\ell^{-3/2}).
\end{equation}
\end{Proposition}

\begin{proof}
(i) Using Cauchy's formula at infinity, we write
\begin{equation}
\label{eq:eqinfty0}
\frac{1}{\ell!}V^{(\ell+2)}(\infty)
=-\frac{1}{2\pi i}\int_{\mathfrak{C}}V''(s)\, s^{\ell-1}\, ds,
\end{equation}
where $\mathfrak{C}$ is the circle $|s|=1+\ell^{-1}$.

Parameterizing $s=(1+\ell^{-1})e^{i\phi}$, we obtain
\begin{equation}
    \label{twostars}
\left|\int_{\mathfrak{C}}V''(s)\, s^{\ell-1}\, ds\right|
\lesssim (1+\ell^{-1})^\ell \int_{-\pi}^{\pi}
\left|V''\!\left((1+\ell^{-1})e^{i\phi}\right)\right|\, d\phi.
\end{equation}

By Theorem~\ref{T:BrillouinReg}(i), formula \eqref{eq:estimv20} with $\epsilon=\ell^{-1}$, we have
\begin{equation}
    \label{onestar}
|V''((1+\ell^{-1})e^{i\phi})|\lesssim \ell^{1/2},
\end{equation}
uniformly in $\phi$. Hence
\[
\int_{-\pi}^{\pi}
\left|V''((1+\ell^{-1})e^{i\phi})\right|\, d\phi
\lesssim \ell^{1/2},
\]
and since $(1+\ell^{-1})^\ell \to e$, we obtain from \eqref{twostars}
\[
\left|\int_{\mathfrak{C}}V''(s)\, s^{\ell-1}\, ds\right|
\lesssim \ell^{1/2}.
\]

It follows that
\[
\frac{1}{\ell!}V^{(\ell+2)}(\infty)=O(\ell^{1/2}),
\]
and integrating twice in $z$ yields
\[
a_\ell=O(\ell^{-3/2}).
\]
(ii) The proof is the same as (i), using now  Theorem \ref{T:BrillouinReg} (ii).

\end{proof}

\subsubsection{Proof of Theorem \ref{P:decay2}} 

\begin{proof} Suppose for contradiction $\lim_{\eps \to 0}\epsilon^{1/2 - \alpha}
V''\bigl(e^{i\phi}(1+\epsilon)\bigr)$ is bounded or converges to $0$ for some $\alpha>0$. Then the factor $\epsilon^{1/2}$ can be replaced by $\epsilon^{1/2-\alpha}$ in in \eqref{eq:estimv20} for all potentials satisfying Assumptions~\ref{A:generic}, respectively in \eqref{limVsec} for all potentials satisfying Assumption \ref{A:generic4}. This implies, by reworking the proof of Proposition \ref{P:decay0} (i), respectively (ii), that $a_\ell=O\left(\ell^{-3/2-\alpha}\right)$, contradicting  Theorem \ref{T:opt1} {\bf A}(ii), respectively $a_\ell=o\left(\ell^{-3/2-\alpha}\right)$ contradicting Theorem \ref{T:opt1} {\bf B}(ii).
\end{proof}

\subsection{Additional decay when there is a ``sharp maximal peak''}
\label{pf2}

In what follows we impose an additional geometric
assumption on the maximal peak of the flat planet $T$
defined in Section~\ref{2dplanet}.
We assume that the maximal peak is unique and isolated, occuring at $\theta=\theta_0$,
as in Section~\ref{noationa}, and that, in a neighborhood
of $e^{i\theta_0}$, the planet is contained in a wedge:
\begin{equation}
\label{eq:wedge}
\mathcal W=
\left\{
e^{i\theta_0}(1-\tau e^{is}):
0\leq\tau<\epsilon,\ |s|\leq\gamma
\right\},
\qquad 0<\gamma<\frac{\pi}{2}.
\end{equation}
We choose $\epsilon>0$ sufficiently small that
\[
T\cap\{|w|>1-\epsilon/2\}\subset\mathcal W.
\]
The symmetrized planet $\mathfrak T$ has a
corresponding conjugate peak at $e^{-i\theta_0}$.

\begin{proposition}\label{P:decay}
Under the assumptions above, we have
$a_\ell=O(\ell^{-5/2})$.
\end{proposition}

{\bf Note.} Under an appropriate one-sided Lebesgue
condition at the maximal peak, the estimate can
similarly be improved to $o(\ell^{-5/2})$;
we omit the details.

\begin{proof}
It is convenient to take $\phi\in[-\pi,\pi)$,
since $V$ is analytic in the exterior of $\mathfrak T$.

We split the flat planet $T$ into
\[
T_1=T\cap D(0,1-\tfrac{\epsilon}{2}),
\qquad
T_0=T\setminus T_1.
\]
By the isolation and wedge assumptions, $T_0$
is contained in the wedge \eqref{eq:wedge}.
Denote by $V_1$ and $V_0$ the potentials obtained
by integrating over $T_1$ and $T_0$, respectively.
The coefficients of $V_1$ are exponentially small,
since $T_1$ is contained in a disk of radius
strictly less than $1$.
We therefore focus on the coefficients of $V_0$.

Proceeding as in the proof of
Proposition~\ref{P:decay0}, the contour of integration
used to calculate the coefficients of $V_0''$ can be
deformed to surround the wedge and its complex
conjugate at distance $O(\ell^{-1})$.
The contribution from the circular arc is
exponentially small. The ray contributions are
treated similarly (with coincident rays at a polar
peak), so we consider only one.

By rotation symmetry, we may assume that the
corresponding ray $L$ originates at $1+\ell^{-1}$.
We parameterize it by
\[
s(t)=1+\ell^{-1}-te^{-i\gamma},
\qquad 0\leq t\leq t^*,
\]
where $t^*$ is the smallest positive solution of
$|s(t)|=1-\epsilon/2$. A direct calculation gives
\[
t^*=\frac{\epsilon}{2\cos\gamma}
+O(\epsilon^2)+O(\ell^{-1}).
\]
Along this ray we write
\[
J_0=\int_0^{t^*}
V_0''(s(t))s'(t)s(t)^{\ell-1}\,dt.
\]

We estimate $V_0''$ using \eqref{eq:difdv}. The wedge condition
implies that the contour remains at distance
$\gtrsim\ell^{-1}$ from the symmetrized support of
$V_0$. Noting that
\(
2r\sin\theta\leq
|s-re^{i\theta}|+|s-re^{-i\theta}|,
\)
the two terms in \eqref{eq:difdv}, after integration over $T_0$,
are $O(\ell^{1/2})$ and $O(\log\ell)$, respectively.
Integration over the truncated planet therefore
gives, uniformly along $L$,
\[
|V_0''(s)|\lesssim\ell^{1/2}.
\]
This estimate also holds at polar peaks. The parameterization of $L$ gives
\[
|s(t)|^{\ell-1}\lesssim e^{-c\ell t},
\qquad 0\leq t\leq t^*,
\]
for some $c>0$ independent of $\ell$ implying
\[
|J_0|\lesssim
\ell^{1/2}\int_0^{t^*}e^{-c\ell t}\,dt
=O(\ell^{-1/2}),
\]
Two integrations give
\( 
a_\ell=O(\ell^{-5/2}),
\)
as claimed.
\end{proof}
\begin{note}
  Instead of wedge--shaped containment, one may assume smoother geometry of the
envelope of the peak, for instance parabolic profiles. Perhaps surprisingly, in
general this leads to \emph{slower} decay rates. In the absence of additional
regularity of the density, smoother topography at the peak results in a more rapid
shrinking of the radius of analyticity of $V$ as the peak is approached. The worst
decay, $o(\ell^{-\frac32})$, occurs when the peak occupies an entire interval in
$\theta$\footnote{Of sufficient extent so that the mass of the peak region
contributes to the spherical harmonic coefficients at practically attainable
orders and accuracy.}, in which case the $o(\ell^{-\frac52})$ behavior is due to the
non--smoothness of the density.

\end{note}
\section{Proof of Theorem~\ref{T:iff}}
\label{pf4}

\begin{proof}
Let $\kappa\in(\kappa_0,1)$, with $\kappa_0$ defined by
\eqref{kappa0}. Since the SHE \eqref{eq:eqvs} for the gravitational potential converges at $z_c$, it
converges to an analytic function in
\[
\{z:|z|>|z_c|\}.
\]
Denote this function by $\widehat V(z)$. Since the series agrees with
$V(z)$ for $|z|>1$, the function $\widehat V$ is the analytic
continuation of $V$ to $\{z:|z|>|z_c|\}$, and hence in particular to
the annulus $A=\{z:\kappa_0<|z|<1\}$.

Recall $V_\kappa$ and $V^\kappa$ denote the potentials of $P_\kappa$ (the inner planet) and
$P^\kappa$ (the outer planetary annulus), and $V=V_\kappa$ and $V^\kappa$ . The possible complex singularities of $V^\kappa$ lie in the
symmetrized set $\mathfrak T^\kappa$. We first show that if, after analytic
continuation, $\hat{V}^\kappa$ has singularities in $A$, they must lie on $C=
\partial\mathfrak T^\kappa
\cap
\partial\mathfrak T_\kappa$ (the common boundary of the complexified inner planet and outer annular planet).

\begin{lemma}\label{L:L10}
The functions $V_\kappa$ and $V^\kappa$ admit analytic continuations
$\widehat V_\kappa$ and $\widehat V^\kappa$, respectively, to
$A\setminus C$, with continuous boundary values on both sides of $C$.
\end{lemma}

\begin{proof}
The potential $V$ is analytic in $\mathfrak T^c$. Similarly,
$V_\kappa$ is analytic in $(\mathfrak T_\kappa)^c$, while
$V^\kappa$ is analytic in $(\mathfrak T^\kappa)^c$.

Since
\[
V_\kappa(z)=V(z)-V^\kappa(z), \qquad |z|>1,
\]
and $\widehat V$ is analytic in $A$, the right-hand side gives an
analytic continuation of $V_\kappa$ to
\[
A\cap(\mathfrak T^\kappa)^c.
\]
On the other hand, $V_\kappa$ itself is analytic in
\[
A\cap(\mathfrak T_\kappa)^c.
\]
The two analytic functions agree where both coincide with the original
$V_\kappa$, and hence, by analytic continuation, together define an
analytic function $\widehat V_\kappa$ on
\[
\bigl[A\cap(\mathfrak T^\kappa)^c\bigr]
\cup
\bigl[A\cap(\mathfrak T_\kappa)^c\bigr]
=
A\setminus C.
\]

The same argument, starting from
\[
V^\kappa(z)=V(z)-V_\kappa(z), \qquad |z|>1,
\]
gives an analytic continuation $\widehat V^\kappa$ to $A\setminus C$.

Theorem~\ref{T:BrillouinReg}(iv), after translation and rescaling,
shows that these continuations have continuous boundary values on both
sides of $C$.
\end{proof}

We now recover $\widehat V^\kappa$ from its jump across $C$. Let
$\mathcal C$ be a smooth closed curve surrounding
$\mathfrak T^\kappa$, oriented clockwise, and let $z$ lie in its
exterior. Cauchy's formula at infinity gives
\begin{equation}
\label{eq:eqinfty}
\widehat V^\kappa(z)
=
\frac{1}{2\pi i}
\int_{\mathcal C}
\frac{\widehat V^\kappa(s)}{s-z}\,ds.
\end{equation}

By Lemma~\ref{L:L10}, $\widehat V^\kappa$ is analytic in
$A\setminus C$ and has continuous boundary values on both sides of
$C$. We may therefore deform $\mathcal C$ onto the two sides of $C$.
Orient $C$ counterclockwise along the circle $|s|=\kappa$. Taking the
limit in \eqref{eq:eqinfty} gives
\begin{equation}
\label{eq:eqinftyshinked}
\widehat V^\kappa(z)
=
-\frac{1}{2\pi i}
\int_C
\frac{
\widehat V^{\kappa,+}(s)-\widehat V^{\kappa,-}(s)
}{s-z}\,ds.
\end{equation}

Since
\[
\widehat V
=
\widehat V_\kappa+\widehat V^\kappa
\]
is analytic across $C$, it has no jump there. Therefore
\[
\widehat V^{\kappa,+}(s)-\widehat V^{\kappa,-}(s)
=
-\bigl(
\widehat V_\kappa^+(s)-\widehat V_\kappa^-(s)
\bigr)
=
\widehat V_\kappa^-(s)-\widehat V_\kappa^+(s)
=
\Phi(s).
\]
Substituting this into \eqref{eq:eqinftyshinked}, we obtain
\[
\widehat V^\kappa(z)
=
-\frac{1}{2\pi i}
\int_C\frac{\Phi(s)}{s-z}\,ds,
\qquad |z|>1.
\]
Since $V=V_\kappa+V^\kappa$ for $|z|>1$, it follows that
\[
V(z)
=
V_\kappa(z)
-
\frac{1}{2\pi i}
\int_C\frac{\Phi(s)}{s-z}\,ds,
\qquad |z|>1.
\]
This is \eqref{eq:eqV} and proves part~(a). The right-hand side is
analytic for $|z|>\kappa$ and gives the analytic continuation of $V$
to this domain.

For part~(b), the coefficient of $1/z$ in the expansion of $V$ at
infinity is $G$ times the total mass of $P$. By part~(a), $V$ is
completely determined by $V_\kappa$, and therefore so is the total
mass of $P$.

For part~(c), suppose that, for some $\kappa\in(\kappa_0,1)$,
$V$ satisfies \eqref{eq:eqV} and that its right-hand side is analytic
for $|z|>\kappa$. It therefore gives an analytic continuation of $V$
to this domain. The Taylor series of $V$ at infinity consequently
converges for $|z|>\kappa$. Since this series agrees with the SHE
\eqref{eq:eqvs} for $|z|>1$, uniqueness of the Taylor expansion at
infinity shows that it is precisely the SHE of $V$. Hence the SHE
converges for $|z|>\kappa$.
\end{proof}

 \section{The balayage measure on the Brillouin sphere and a formula for the circular averages $\mu$} 
 \subsection{Proof of Theorem~\ref{mufromV}}
\label{pfmufromV}

\begin{proof}
Given the potential $V$ with its expansion \eqref{eq:eqvs}, definition \eqref{tildeV} gives
\begin{equation}
\label{eq:ck}
\widetilde V(\zeta)
=
G^{-1}\sum_{\ell=0}^\infty a_{\ell+1}\zeta^\ell,
\qquad |\zeta|<1.
\end{equation}

Let $\mu \in L^2([-1,1])$ satisfy \eqref{eq:balayage1}. With definition
\eqref{tildeV}, we have
\begin{equation}
\label{tildeV2}
\widetilde V(\zeta)
=
\int_{-1}^1
\frac{\mu(x)\,dx}{\sqrt{1-2\zeta x+\zeta^2}},
\qquad |\zeta|<1.
\end{equation}

Expanding the kernel in \eqref{tildeV2} by the generating function for
the Legendre polynomials and applying dominated convergence, we obtain
\[
\widetilde V(\zeta)
=
\sum_{\ell=0}^\infty
\left(
\int_{-1}^1\mu(x)P_\ell(x)\,dx
\right)\zeta^\ell.
\]
Comparison with \eqref{eq:ck} gives
\begin{equation}
\label{alPl}
G^{-1}a_{\ell+1}
=
\int_{-1}^1\mu(x)P_\ell(x)\,dx,
\qquad \ell\geq0.
\end{equation}

Since $\mu\in L^2([-1,1])$ (see Note \ref{sevenp1}(c)), it
has the Legendre expansion
\begin{equation}
\label{sermu}
\mu(x)
=
\sum_{\ell=0}^\infty d_\ell P_\ell(x)
\end{equation}
in $L^2([-1,1])$; see Proposition~\ref{P:legendre}(ii). By
\eqref{alPl} and the orthogonality relation \eqref{eq:ortho},
\[
G^{-1}a_{\ell+1}
=
d_\ell\int_{-1}^1P_\ell^2(x)\,dx
=
d_\ell\frac{2}{2\ell+1}.
\]
Hence
\begin{equation}
\label{eq:eqdk}
d_\ell
=
G^{-1}\left(\ell+\frac12\right)a_{\ell+1}.
\end{equation}
With \eqref{sermu}, we see that we can represent $\mu \in L^2([-1,1])$ in terms of the coefficients of $V$'s expansion; precisely, $$\mu(x) = \sum_{\ell=0}^\infty G^{-1} (\ell + \frac{1}{2}) a_{\ell+1} P_{\ell}(x).$$
However, more can be derived. 

For $0<\rho<1$, define the Abel means of the Legendre expansion by
\begin{equation}
\label{eq:expmu1}
\mu_\rho(x)
:=
\sum_{\ell=0}^\infty d_\ell\rho^\ell P_\ell(x).
\end{equation}
Then
\begin{equation}
\label{eq:AbelL2}
\mu_\rho\longrightarrow\mu
\qquad\text{in }L^2([-1,1])
\quad\text{as }\rho\uparrow1.
\end{equation}
Indeed, by orthogonality,
\[
\|\mu-\mu_\rho\|_{L^2}^2
=
\sum_{\ell=0}^\infty
|d_\ell|^2(1-\rho^\ell)^2\|P_\ell\|_{L^2}^2,
\]
and the right-hand side tends to zero by dominated convergence.

For $x=\cos\theta$, with $0<\theta<\pi$, we use the Stieltjes
representation of the Legendre polynomials
(see \cite{Szeg\"o}, (4.8.17), and also \cite{Maric}):
\begin{equation}
\label{Stieltjes}
P_\ell(\cos\theta)
=
\frac{2}{\pi}
\Im\left\{
e^{i(\ell+1)\theta}
\int_0^1
\frac{t^\ell\,dt}
{\sqrt{1-t}\sqrt{1-te^{2i\theta}}}
\right\}.
\end{equation}
Substituting \eqref{eq:eqdk} and \eqref{Stieltjes} into
\eqref{eq:expmu1}, we obtain
\begin{align}
\mu_\rho(\cos\theta)
&=
\frac{2}{\pi}
\Im\left[
e^{i\theta}
\sum_{\ell=0}^\infty
G^{-1}\left(\ell+\frac12\right)a_{\ell+1}
\int_0^1
\frac{(\rho e^{i\theta}t)^\ell}
{\sqrt{1-t}\sqrt{1-te^{2i\theta}}}
\,dt
\right].
\label{interfrom}
\end{align}

Define
\begin{equation}
\label{eq:linkVH}
H(\zeta)
=
G^{-1}
\sum_{\ell=0}^\infty
\left(\ell+\frac12\right)a_{\ell+1}\zeta^\ell,
\qquad |\zeta|<1.
\end{equation}
By \eqref{eq:ck}, direct differentiation gives
\[
H(\zeta)
=
\zeta\widetilde V'(\zeta)
+\frac12\widetilde V(\zeta),
\]
which is \eqref{formH}.

For $0<\rho<1$, the series in \eqref{interfrom} may be interchanged
with the integral, giving
\begin{equation}
\label{formularho}
\mu_\rho(\cos\theta)
=
\frac{2}{\pi}
\Im\left[
e^{i\theta}
\int_0^1
\frac{H(\rho t e^{i\theta})}
{\sqrt{1-t}\sqrt{1-te^{2i\theta}}}
\,dt
\right],
\qquad 0<\theta<\pi.
\end{equation}

By Theorem~\ref{T:BrillouinReg}(iv), $V$ and $V'$ extend continuously
to the Brillouin sphere. Under the change of variables $\zeta=1/z$,
in \eqref{formH}, we see that $H$ extends continuously to
$\overline{\mathbb D}$. Hence, for each fixed
$\theta\in(0,\pi)$, dominated convergence in \eqref{formularho} gives
\[
\lim_{\rho\uparrow1}\mu_\rho(\cos\theta)
=
\frac{2}{\pi}
\Im\left[
e^{i\theta}
\int_0^1
\frac{H(te^{i\theta})}
{\sqrt{1-t}\sqrt{1-te^{2i\theta}}}
\,dt
\right]
=
\mu_*(\cos\theta),
\]
as defined in \eqref{muofcost}. Thus $\mu_\rho(x)\to\mu_*(x)$ for every $x\in(-1,1)$.

It remains to identify the endpoint values. Since
\[
P_\ell(1)=1,
\qquad
P_\ell(-1)=(-1)^\ell,
\]
we have, by \eqref{eq:expmu1} and \eqref{eq:linkVH},
\[
\mu_\rho(1)=H(\rho),
\qquad
\mu_\rho(-1)=H(-\rho).
\]
The continuity of $H$ on $\overline{\mathbb D}$ therefore gives
\[
\lim_{\rho\uparrow1}\mu_\rho(1)=H(1),
\qquad
\lim_{\rho\uparrow1}\mu_\rho(-1)=H(-1).
\]
This gives the endpoint values in \eqref{muofcostendpoints}.

We next verify continuity at the endpoints. Although the integral in
\eqref{Stieltjes} has a logarithmic singularity as $\theta\downarrow0$,
its imaginary part has a finite limit. Indeed, writing
$t^\ell=1+(t^\ell-1)$, the second term has an integrable limit, while
for the first term
\[
e^{i\theta}
\int_0^1
\frac{dt}
{\sqrt{1-t}\sqrt{1-te^{2i\theta}}}
=
\log\cot\frac{\theta}{2}+\frac{i\pi}{2}.
\]
Consequently,
\[
\lim_{\theta\downarrow0}
\frac{2}{\pi}
\Im\left\{
e^{i(\ell+1)\theta}
\int_0^1
\frac{t^\ell\,dt}
{\sqrt{1-t}\sqrt{1-te^{2i\theta}}}
\right\}
=
1
=
P_\ell(1).
\]
The corresponding statement at $\theta=\pi$ follows by symmetry.
Together with the continuous extension of $H$ to
$\overline{\mathbb D}$, this shows that the function $\mu_*$ defined
by \eqref{muofcost} and \eqref{muofcostendpoints} is continuous on
$[-1,1]$.

Finally, by \eqref{eq:AbelL2},
\[
\mu_\rho\longrightarrow\mu
\qquad\text{in }L^2([-1,1]),
\]
whereas we have shown that
\[
\mu_\rho(x)\longrightarrow\mu_*(x)
\qquad\text{for every }x\in[-1,1].
\]
Every $L^2$-convergent sequence has a subsequence converging pointwise almost
everywhere to its $L^2$ limit. Applying this to a sequence
$\rho_n\uparrow1$ and using the pointwise limit above yields
\[
\mu_*=\mu
\qquad\text{a.e. on }[-1,1].
\]
Thus, \eqref{eq:balayage1} is satisfied with $\mu := \mu_*$, proving the theorem.
\end{proof}

\subsection{Proof of Theorem\, \ref{T:sigmaU}: the balayage measure in terms of the potential}\label{pfsigmau}

We assumed $\rho\in L^\infty (P)$. By allowing $\rho$ to vanish outside $P$ we can assume $P$ is the unit ball (recall that we assumed the Brillouin radius is $1$). 

We denote by \(\{Y_{lm}\}_{l\ge 0, |m|\le l}\) the orthonormal set of spherical harmonic functions $Y_{lm}$. If $\mathbf z$ has spherical coordinates $\rho,\theta,\lambda$ we denote functions $f(\mathbf z)$ also by $f(\rho,\theta,\phi)$, and $Y_{lm}\left(\frac{\mathbf z}{|\mathbf z|}\right)$ will be also denoted by $Y_{lm}(\theta,\lambda)$.

Consider the potential for an arbitrary point $\mathbf z$ outside the planet \eqref{formU}.

Using the expansion of the inverse distance in the integrand in \eqref{formU} (see \S\ref{expinvdidt}, and in particular \eqref{recdistexpansion} we obtain for $U$ the spherical harmonic expansion
\begin{equation}
\label{epandV}
U(r,\theta ,\phi) =G\sum _{l=0}^{\infty }\frac 1{r^{l+1}}\, U_l(\theta,\phi)\ \ \text{where }U_l(\theta,\phi)=\sum _{m=-l}^{l}{I_{lm}}Y_{l,m}(\theta ,\phi )
\end{equation}
where $(r,\theta ,\phi)$ are the spherical coordinates of $\mathbf z$ and $I_{lm}$ are given by
\begin{equation}
\label{Ilm}
 I_{lm}=\frac{4\pi }{2l+1}\int_P \rho (\mathbf y)\, |\mathbf y|^{l}\, \ \overline{Y_{l,m}}\left(\frac{\mathbf y}{|\mathbf y|}\right) \, d\mathbf y
 \end{equation}

 It is known that the series \eqref{epandV} converges for $r>1$. Therefore, $U(r,\theta ,\phi)$ is an analytic function for $r\in \CC$ with $|r|>1$. In fact, more is known: after a rotation of the coordinate system, 
Proposition\,\ref{P:decay0} applies, and, under Assumption \ref{A:generic} we have 
$|U_{\ell}(\theta,\phi)|=O(\ell^{-3/2})$.
 
Consider the associated balayage measure $\sigma$ given by \eqref{Ubalayage1comp}. Since we assumed $\rho\in L^\infty(P)$, then we also have $\sigma\in L^\infty(S^2)$ (see \S\ref{sigmaisLinfty}).

{\em Proof.}

The expansion \eqref{epandV} implies that
\begin{equation}
    \label{serTU}
\tilde{U}(\zeta,\theta,\lambda)=\sum _{l=0}^{\infty }\zeta^{l}\, U_l(\theta,\lambda)
\end{equation}
and the series converges for $\zeta\in\CC$ with $|\zeta|<1$.

Since $\sigma\in L^\infty(S^2)$ it follows that $\sigma\in L^2(S^2)$ and it has an $L^2$ expansion in spherical harmonics
\begin{equation}
    \label{sersigma}
\sigma(\mathbf u)=\sum_{l=0}^\infty\sum_{m=-l}^l\sigma_{lm}Y_{lm}(\mathbf u)
\end{equation}
which we use in \eqref{Ubalayage1comp}, together with \eqref{recdistexpansion}, and obtain
\[
U(\mathbf z)=G\int_{S^2} \left[\sum_{l,m}\sigma_{lm}Y_{lm}(\mathbf u)\right]\, \left[ \sum _{l,m}\frac{4\pi }{2l+1}\frac1{r^{l+1}}\, Y_{l,m}\left(\frac{\mathbf z}{|\mathbf z|}\right)\overline{Y_{l,m}}\left(\mathbf u\right)\right]\, dA(\mathbf u)
\]
\[
=G\sum_{l,m}\frac{4\pi }{2l+1}\frac1{r^{l+1}}\sigma_{lm}Y_{l,m}\left(\frac{\mathbf z}{|\mathbf z|}\right)
\]
which compared to \eqref{epandV} yields that
\[\frac{4\pi }{2l+1}\sigma_{lm}=I_{lm}
\]
therefore
\begin{equation}
    \label{sigmaI}
2\pi\, \sigma_{lm}=\left(l+\frac12\right)I_{lm}
\end{equation}

Let $\eta\in[0,1)$. Since \eqref{sersigma} converges in $L^2$ it follows that the series
\begin{equation}
    \label{sersigmaeta}
\sigma(\mathbf u;\eta)=\sum_{l=0}^\infty\eta^l\sum_{m=-l}^l\sigma_{lm}Y_{lm}(\mathbf u)
\end{equation}
is convergent, $\sigma(\mathbf u;\eta)$ analytic for $|\eta|<1$, and
\( \sigma(\mathbf u)=\lim_{\eta \nearrow  1} \sigma(\mathbf u;\eta)\), a limit which exists in $L^2(S^2)$ since by Parseval and dominated convergence,
$$ \|\sigma(\cdot;\eta)-\sigma\|_{L^2}^2 = \sum_{\ell,m} |1-\eta^\ell|^2|\sigma_{\ell m}|^2 \to0. $$
 Using \eqref{sigmaI}, then \eqref{epandV} and \eqref{serTU} we find that
\begin{equation}
    \label{eq20}
2\pi\, \sigma(\mathbf u;\eta)=\sum_{l=0}^\infty\eta^l\left(l+\frac12\right)\sum_{m=-l}^lI_{lm}Y_{lm}(\mathbf u)=\sum_{l=0}^\infty\eta^l\left(l+\frac12\right)U_l(\mathbf u)=\eta\partial_\eta\tilde{U}(\eta,\mathbf u)+\frac12 \tilde{U}(\eta,\mathbf u)
\end{equation}

Theorem\,\ref{T:BrillouinReg}(iv) shows that if the planet satisfies Assumptions~\ref{A:generic}, then $V(z):=U(z \mathbf{e_3})$, where $\mathbf{e_3} = (0,0,1)$  and $V'(z)=\nabla U \cdot \mathbf{e_3}$ are continuous in the closed complement of the unit disk, more precisely $\lim_{z \to z_0,|z|>1} V(z)$ and $\lim_{z \to z_0,|z|>1} V'(z)$ exist for all $|z_0|\ge 1$. Rotating the coordinate system, the same is true for $V(z):=U(z\mathbf u)$ for any unit vector $\mathbf u$. (Indeed, if $R$ is a rotation so that $\mathbf u =R\mathbf e_3$, then $U(\mathbf z)$ equals the potential $U_R(z\mathbf e_3)$ associated to the rotated density $\rho\circ R$). Therefore, the expression on the right hand side of \eqref{eq20} is continuous in $|\eta|\le 1$ and has pointwise limit as $\eta \nearrow  1 $. 

It follows from \eqref{eq20} that $\sigma(u;\eta)$ has a pointwise
limit as $\eta\uparrow1$ for every $u\in S^2$. Since
$\sigma(\,\cdot\,;\eta)\to\sigma$ in $L^2(S^2)$, the result follows.
\

\

\section{Appendix}

\subsection{Legendre polynomials}

\subsubsection{Some classical results on Legendre polynomials}

This section summarizes some classical results on Legendre polynomials.

\begin{Proposition}[\cite{Wiki1}, \cite{Szeg\"o}, (7.21.1), p.~164]
\label{P:legendre}
{\rm
(i) The Legendre polynomials $\{P_n(x)\}_{n \in \mathbb{N}\cup \{0\}}$ admit the generating function
\begin{equation}
  \label{defG}
  G(\zeta;x)= \frac{1}{\sqrt{1 - 2x\zeta + \zeta^{2}}},
  \end{equation}
as
\begin{equation}
  \label{eq:gener}
  \frac{1}{\sqrt{1 - 2x\zeta + \zeta^{2}}}
  = \sum_{\ell =0}^{\infty} P_{\ell }(x)\,\zeta^{\ell },
\end{equation}
for $|x| \leq 1$ and $|\zeta|<1$.

(ii) $\{P_n\}_{n\ge 0}$ form a complete orthogonal system in $L^2([-1,1])$, satisfying
\begin{equation}
  \label{eq:ortho}
  \int_{-1}^{1} P_m(x) P_n(x)\,dx 
  \begin{cases}
    0, & n \neq m,\\[4pt]
    \displaystyle \frac{2}{2n+1}, & n = m.
  \end{cases}
\end{equation}
(iii) Moreover,
\begin{equation}
  \label{eq:intrep}
  \sup_{\ell \ge 0,\, x \in [-1,1]} |P_\ell(x)| \le 1.
\end{equation}
}
\end{Proposition}

The Bernstein inequality for Legendre polynomials is (\cite{Szeg\"o} Theorem 7.3.3, \cite{Bernstein}, \cite{Maric}) :

\begin{equation}
\label{eq:szegoineq}
 \sqrt{\sin\theta}\left| P_n\left(\cos\theta\right)\right|\le \sqrt{\frac{2}{n\pi}},\ \ \ \ \ \text{for }\theta\in[0,\pi] 
 \end{equation}

 We also need Hilb's asymptotic formula (see \cite{Szeg\"o}, also \cite{nist} 18.15.6)
 \begin{equation}
     \label{HilbShort}
     P_n(\cos\theta)=\sqrt{\frac{\theta}{\sin\theta}}\, J_0\left(\left(n+\frac12\right)\theta\right)+O\left(n^{-3/2}\right)\, \min\{\theta^{1/2},(\pi-\theta)^{1/2}\}
 \end{equation}
with uniform estimates for errors for $cn^{-1}\le\theta\le \pi-cn^{-1}$.

Using the asymptotic formula for Bessel functions in \eqref{HilbShort} we get 
 \begin{equation}
     \label{Hilb0}
     P_n(\cos\theta)=\sqrt{\frac{2}{n\pi\sin\theta}}\ \cos\left[\left(n+\frac12\right)\theta-\frac\pi 4\right]+o\left(n^{-1/2}\theta^{-1/2}\right) \ \ \text{ for }n\theta\gg 1\ \text{ and }\theta<\pi-\epsilon
 \end{equation}
  (see also \cite{Szeg\"o}, Theorem 8.21.2). By symmetry, since $P_n(\cos(\pi-\theta))=(-1)^nP_n(\cos(\theta))$ we obtain
  \begin{equation}
     \label{HilbPi}
     P_n(\cos\theta)=\sqrt{\frac{2}{n\pi\sin\theta}}\ \cos\left[\left(n+\frac12\right)\theta-\frac\pi 4\right]+o\left(n^{-1/2}(\pi-\theta)^{-1/2}\right) \ \ \text{ for }n(\pi-\theta)\gg 1\ \text{ and }\theta>\epsilon
 \end{equation}
 Relations \eqref{Hilb0} and \eqref{HilbPi} can be written in a unified way as 
 \begin{equation}
\label{eq:Piscos}
P_{\ell }(\cos \theta )=\sqrt {\frac {2}{\pi \ell \sin \left(\theta \right)}}\cos \left[\left(\ell +{\tfrac {1}{2}}\right)\theta -{\tfrac {\pi }{4}}\right]+{\mathcal {O}}\left(\ell ^{-3/2}\right);\ \ \text{ for }\ell\theta\gg 1
\end{equation}

\subsubsection{Estimates of Legendre polynomials in a complex neighborhood of $[-1,1]$}\label{Sec2p1}

Consider the generating function \eqref{defG} for the Legendre polynomials. 
For $x\in[-1,1]$,  with the natural branch of the square root which is positive when $\zeta$ and $x$ are small and positive, $G(\cdot;x)$ 
extends to a holomorphic map in the open unit disk $\mathbb{D}$.

  Let $C$ be a positively oriented, simple closed contour around 0 inside $\mathbb{D}$. From \eqref{eq:gener}, using Cauchy's formula, we have
\begin{equation}
  \label{eq:Cauchy}
  P_n(x)
  = \frac{1}{2\pi i}
    \oint_C \frac{ds}{s^{n+1}\sqrt{s^2 - 2xs + 1}},
\end{equation}

\begin{Proposition}\label{estimPn} 
{\rm 
 Let  $x=\cos(\theta+iy)$ where $\theta,y\in\RR$ and $|y|<\ln L$.  For any $y'>|y|$ we have
$$|P_n(\cos(\theta+iy))|\le  K\,  e^{ny'},\ \ \ \ \text{where }K=K_{y,y'}=\left(e^{-y}- e^{-y'}\right)^{-1}$$
}
\end{Proposition}

\begin{proof}
  In \eqref{eq:Cauchy} we choose $y'>|y|$ and $C$ to be the circle  $\{s:|s|=e^{-y'}\}$. Then the roots $s_{\pm}=e^{\pm i(\theta+iy)}$ of $s^2-2 x s+1$ are outside $C$. Furthermore, along $C$, we have $|s^{-n-1}|=e^{(n+1)y'}$ and
  $$\left|s^2-2 x\cos(\theta+iy)+1\right|=\left|s-e^{i(\theta+i y)}\right|\left|s-e^{-i(\theta+i y)}\right|\ge (e^{-y}-e^{-|y|})^2$$
  which used in \eqref{eq:Cauchy} yields 
$$|P_n(\cos(\theta+iy))|\le \frac{e^{(n+1)y'}}{2\pi(e^{-y}-e^{-|y|})}2\pi e^{-y'}=Ke^{ny'}$$
  \end{proof}

\subsection{$r_m(\theta)$ is upper semicontinuous}\label{semic}

\begin{lemma}
Let $P\subset\mathbb{R}^3$ be compact, and define
\[
 r_m(\theta)
 :=
 \max\left\{
 r\geq 0:
 r\,\mathbf{u}(\theta,\lambda)\in P
 \text{ for some }\lambda\in[0,2\pi]
 \right\},
\]
where
\[
 \mathbf{u}(\theta,\lambda)
 =
 (\sin\theta\cos\lambda,\,
  \sin\theta\sin\lambda,\,
  \cos\theta).
\]
Then $r_m$ is upper semicontinuous on $[0,\pi]$, and hence
Borel measurable.
\end{lemma}

\begin{proof}
Let $\theta_n\to\theta$. We show that
\[
 \limsup_{n\to\infty} r_m(\theta_n)
 \leq r_m(\theta).
\]
Passing to a subsequence, if necessary, we may assume that
\[
 r_m(\theta_n)\longrightarrow
 L:=\limsup_{n\to\infty}r_m(\theta_n).
\]
For each $n$, choose $\lambda_n\in[0,2\pi]$ such that
\[
 r_m(\theta_n)\mathbf{u}(\theta_n,\lambda_n)\in P.
\]
Since $[0,2\pi]$ is compact, after passing to a further subsequence
we may assume that
\[
 \lambda_n\longrightarrow\lambda
\]
for some $\lambda\in[0,2\pi]$. By continuity of $\mathbf{u}$,
\[
 r_m(\theta_n)\mathbf{u}(\theta_n,\lambda_n)
 \longrightarrow
 L\mathbf{u}(\theta,\lambda).
\]
Since $P$ is closed, it follows that
\[
 L\mathbf{u}(\theta,\lambda)\in P.
\]
Therefore, by the definition of $r_m(\theta)$,
\[
 L\leq r_m(\theta).
\]
Thus
\[
 \limsup_{n\to\infty}r_m(\theta_n)
 \leq r_m(\theta),
\]
so $r_m$ is upper semicontinuous. In particular, $r_m$ is Borel
measurable.
\end{proof}

\subsection{Expansion of the inverse distance in terms of spherical harmonics}\label{expinvdidt} 

Let $\mathbf r,\mathbf r^\prime \in\RR^3$. Assume $r:=|\mathbf r|>|\mathbf r^\prime|:=r^\prime$. Then, using \( \mathbf{r}\cdot\mathbf{r}^{\prime }=r r'\cos\psi\) and the generating function for the Legendre polynomials \eqref{eq:gener} reciprocal distance function can be written as
 \begin{equation}
     \label{eq2}
 \frac{1}{|\mathbf{r}-\mathbf{r}^{\prime }|}=\frac{1}{r}\left[1-2\left(\frac{r^{\prime }}{r}\right)\cos \psi +\left(\frac{r^{\prime }}{r}\right)^{2}\right]^{-1/2}=\sum _{l=0}^{\infty }\frac{(r^{\prime })^{l}}{r^{l+1}}P_{l}(\cos \psi )
 \end{equation}
The Spherical Harmonic Addition Theorem reads (see \cite{nist} 14.30.E9)
  \[
 P_{l}(\cos \psi )=\frac{4\pi }{2l+1}\sum _{m=-l}^{l}Y_{l,m}(\theta ,\varphi )\overline{Y_{l,m}}(\theta ^{\prime },\varphi ^{\prime })
 \]
where $\cos \psi=cos\theta\cos\theta'+\sin\theta\sin\theta'
\cos\left(\phi-\phi'\right)$ is the angle between the vectors of spherical coordinates \((1,\theta ,\varphi )\) and  \((1,\theta' ,\varphi' )\).
Using the addition theorem in \eqref{eq2} we obtain
\begin{equation}
     \label{recdistexpansion}
 \frac{1}{|\mathbf{r}-\mathbf{r}^{\prime }|}=\sum _{l=0}^{\infty }\frac{4\pi }{2l+1}\frac{(r^{\prime })^{l}}{r^{l+1}}\, \sum _{m=-l}^{l}Y_{l,m}\left(\frac{\mathbf r}r\right)\overline{Y_{l,m}}\left(\frac{\mathbf r'}{r'}\right)
 \end{equation}   

 \subsection{The balayage measure $\sigma$ is in $L^\infty(S^2)$}\label{sigmaisLinfty}

 Assume $\rho\in L^\infty(B_1)$ and let $\sigma$ be the balayage measure of the potential associated to the density $\rho$. This measure is obtained using the Poisson kernel:
 \(\sigma(\mathbf u)=\int_{B_1}P(\mathbf y,\mathbf u)d\mathbf y\) which implies that \( |\sigma(\mathbf u)|\le \|\rho\|_\infty \int_{B_1}P(\mathbf y,\mathbf u)d\mathbf y\). Now since \(P(\mathbf y,\mathbf u)=\frac{1-|\mathbf y|^2}{|\mathbf u-\mathbf y|^3}\), its integral can be explicitly calculated (equals $4\pi/3$) hence $\sigma\in L^\infty$.
 
\section*{Statements and Declarations}

\paragraph{Funding}
The work of O.C. was partially supported by the
National Science Foundation under grant DMS-2206241.

\paragraph{Conflict of interest}
The authors have no competing interests to declare
that are relevant to the content of this article.

\paragraph{Data availability}
Data sharing is not applicable to this article as no
datasets were generated or analyzed during the current study.

\paragraph{AI usage disclosure}
O.C. and R.D.C. used ChatGPT Instant for editorial
assistance, assistance in identifying potentially
relevant literature, and independent checks of
algebraic calculations.
No mathematical ideas, arguments, proofs, or results
were solicited from or supplied by AI.

\end{document}